\pdfoutput=1
\documentclass{IEEEtran}

\usepackage{empheq,amssymb}
\usepackage{braket}
\usepackage{dsfont,mathrsfs}
\usepackage{enumerate}
\usepackage{graphicx}
\usepackage{amsmath}
\usepackage{float}
\usepackage{bbm}

\newtheorem{definition} {Definition}

\newtheorem{theorem}    {Theorem}
\newtheorem{lemma}      {Lemma}

\newcommand\EXP[1]{\mathop{\kern0pt \mathds E}{\Set{#1}}}
\newcommand\PR [1]{\mathop{\kern0pt \Pr}{\Set{#1}}}

\newcommand{\leftexp}[2]%
  {\mathop{}%
   \mathopen{\vphantom{#2}}^{#1}%
   \kern-\scriptspace%
   #2}

\makeatletter
\newcommand\SEQ{\@ifstar\SEQB\SEQA}
\makeatother

\newcommand\SEQA[2][T]{\{#2_t$, $t=1,\dots,#1\}}
\newcommand\SEQB[1]{\{#1_1$, $t=1,\dots\}}

\begin{document}
\title {A Sequential Problem in Decentralized Detection with Communication}
\author{Ashutosh~Nayyar,~\IEEEmembership{Student Member,~IEEE,}
        and~Demosthenis~Teneketzis,~\IEEEmembership{Fellow,~IEEE}
\thanks{A. Nayyar is with the Department
of Electrical Engineering and Computer Science, University of Michigan, Ann Arbor,
MI, 48109 USA e-mail: (anayyar@umich.edu).}
\thanks{D. Teneketzis is with the Department
of Electrical Engineering and Computer Science, University of Michigan, Ann Arbor,
MI, 48109 USA e-mail: (teneket@eecs.umich.edu).}}
\maketitle

\begin{abstract}
    A sequential problem in decentralized detection is considered. Two observers can make repeated noisy observations of a binary hypothesis on the state of the environment. At any time, observer 1 can stop and send a final binary message to observer 2 or it may continue to take more measurements. Every time observer 1 postpones its final message to observer 2, it incurs a penalty. Observer 2's operation under two different scenarios is explored. In the first scenario, observer 2 waits to receive the final message from observer 1 and then starts taking measurements of its own. It is then faced with a stopping problem on whether to stop and declare a decision on the hypothesis or to continue taking measurements. In the second scenario, observer 2 starts taking measurements from the beginning. It is then  faced with a different stopping problem. At any time, observer 2 can decide whether to stop and declare a decision on the hypothesis or to continue to take more measurements and wait for observer 1 to send its final message. Parametric characterization of optimal policies  for the two observers are obtained under both scenarios. A sequential methodology for finding the optimal policies is presented. The parametric characterizations are then extended to problem with increased communication alphabet for the final message from observer 1 to observer 2; and to the case of multiple peripheral sensors that each send a single final message to a coordinating sensor who makes the final decision on the hypothesis.
\end{abstract}
\section{Introduction}
     Decentralized detection problems are motivated by applications in large scale decentralized systems such as sensor networks and surveillance networks. In such networks, sensors receive different information about the environment but share a common objective, for example to detect the presence of a target in a surveillance area. Sensors may be allowed to communicate but they are constrained to exchange only a limited amount of information because of energy constraints, data storage and data processing constraints, communication constraints etc. 
     \par
    Decentralized detection problems may be static or sequential. In static problems, sensors make a fixed number of observations about a hypothesis on the state of the environment which is modeled as a random variable \(H\). Sensors may transmit a single message (a quantized version of their observations) to a fusion center which makes a final decision on \(H\). Such problems have been extensively studied since their initial formulation in \cite{Tenney_Detection} (See the surveys in \cite{Tsitsiklis_survey}, \cite{Varshney} and references therein).  In most such formulations, it has been shown that person-by-person optimal decision rules (as defined in \cite{Radner_team}) for a binary hypothesis detection problem are characterized by thresholds on the likelihood ratio (or equivalently on the posterior belief on the hypothesis). Under certain conditions such as large number of identical sensors, it has been shown that it is optimal to use identical quantization rule at all sensors (\cite{Tsitsiklis_large}, \cite{Chamberland04}). A related information-theoretic formulation with constraints on communication from a sensor to a fusion center appears in \cite{Ahlswede_hypothesis}.
    \par
 
 In sequential problems, the number of observations taken by the sensors is
not fixed a priori. Two distinct formulations have been considered
for sequential problems. In one formulation, at each time instant local/peripheral sensors
send a message about their observations to a fusion center/coordinator. At each time instant, the fusion center
decides whether to receive more messages or to declare a
decision on the hypothesis. Thus the fusion center is faced with an optimal stopping problem whereas the peripheral sensors are not faced with an optimal stopping problem. The case where peripheral sensors can only use
their current observation and possibly all past transmissions of all sensors to decide what message to send to the fusion center  has been
studied in \cite{VBP_Detection}. No positive results have been
found in the case when sensors remember their past observations as
well.
\par
A second formulation may be motivated by situations where
       continuous communication with a fusion center is too costly because of the various constraints mentioned
       earlier. In this formulation, each sensor locally decides when to stop taking more measurements and only sends a final message to a fusion center. Each sensor pays a penalty for delaying its final decision.
      The fusion center has to wait to receive the final messages from all sensors and then combine them to produce a final decision on the
      hypothesis. Thus, in this formulation, each local/peripheral sensor is faced with an optimal stopping problem but the coordinator does not have a stopping problem. A version of this problem (called the decentralized Wald problem) was formulated in
      \cite{Dec_Wald} and it was shown that at each time instant, optimal policies for the
      peripheral sensors are described by two thresholds. The computation of these thresholds requires solution of two coupled sets of dynamic programming equations. Similar results were obtained in a continuous time setting in \cite{LaVigna_86}. Although
      this formulation reduces the communication requirements, the
      final decision at the fusion center is made only when all
      sensors have sent their messages. In a similar formulation, the problem of quickest detection of the change of state of a Markov chain was considered in \cite{Teneket_Quickest_Det}.
      \par
      In the problem we consider in this paper, the peripheral sensors as well as a coordinating sensor are faced with optimal stopping problems. The peripheral sensors decide locally when they want to stop taking measurements and send a final message
       to a special coordinating sensor, say S0. The coordinating sensor S0 is faced with a  stopping problem of its own. At any time,
       the coordinating sensor S0 uses
        its own measurements
        and the messages it has received so far to make a decision on whether to stop and declare a final decision on the hypothesis or
         continue to take more measurements and wait for messages from other sensors that have not yet sent a final message.
         As in \cite{Dec_Wald}, each sensor (peripheral sensors and the coordinating sensor) incurs a penalty
         for delaying its final message/decision, and a cost depending on
         S0's final decision on the hypothesis and the true value
         of the hypothesis is incurred in the end. 
\par
 We first consider a simple two sensor version of this problem and obtain a parametric characterization of optimal policies. We prove that at each time instant, an optimal policy of the peripheral sensor is characterized by at most 4 thresholds on its posterior belief on the hypothesis; an optimal policy of the coordinating sensor is characterized by 2 thresholds (on its own posterior belief) that depend on the messages received from the peripheral sensor. This characterization differs from the classical two threshold characterization found in the centralized and the decentralized Wald problems (\cite{Wald}, \cite{Dec_Wald}). The computation of these optimal thresholds is a  difficult problem. We present a sequential methodology that decomposes the overall optimization problem into several smaller problems that may be solved to determine the optimal thresholds at each time instant. We extend our results to a problem with multiple peripheral sensors that send their final message to the coordinating sensor who makes the final decision on the hypothesis. We show that qualitative properties of the optimal policies of the peripheral sensors and the coordinating sensor are same as in the two sensor problem.
         \par
         The rest of the paper is organized as follows. In Section~\ref{sec:PF}, we formulate two versions of our problem with two observers. We obtain qualitative results on the nature of optimal policies for the two sensors in Sections~\ref{sec:QP1} and \ref{sec:QP2}. We present a sequential methodology for computing optimal policies in Section~\ref{sec:Optimal_Policies}. In Section~\ref{sec:Inf_Horizon}, we extend our qualitative results to infinite horizon analogues of our problem. A generalization to more than binary communication alphabet is presented in Section~\ref{sec:enhanced_alpha}. We extend our results to a multiple sensor (more than 2) problem in section \ref{sec:MS}. Finally, we conclude in Section~\ref{sec:con}.   

\emph{Notation:} Throughout this paper, \(X_{1:t}\) refers to the sequence \(X_1, X_2,..,X_t\). Subscripts are used as time index and the superscripts are used as the index of the sensor. We use capital letters to denote random variable and the corresponding lower case letters for their realizations. 


\section{Problem formulation} \label{sec:PF}
\subsection{ The Model}
   Consider a binary hypothesis problem where the true hypothesis is modeled as a random variable \(H\) taking values 0 or 1 with known prior probabilities:
   \[P(H=0) = p_0; \hspace{10pt} P(H=1) = 1-p_0\]
 Consider two observers: Observer 1 (O1) and Observer 2 (O2). We assume that each observer can make noisy observations of the true hypothesis. Conditioned on the hypothesis \(H\), the following statements are assumed to be true: \\
 1. The observation of the \(i^{th}\) observer at time \(t\), \((Y^i_t)\) (taking values in the set \(\mathcal{Y}^i\)), either has a discrete distribution \((P^i_t(.|H))\) or admits a probability density function \((f^i_t(.|H))\). \\
2. Observations of the \(i^{th}\) observer at different time instants are conditionally independent given \(H\).\\
3. The observation sequences at the two observers are conditionally independent given \(H\). \\
\begin{figure}[H]
\begin{center}

\includegraphics[height=6cm,width=8cm]{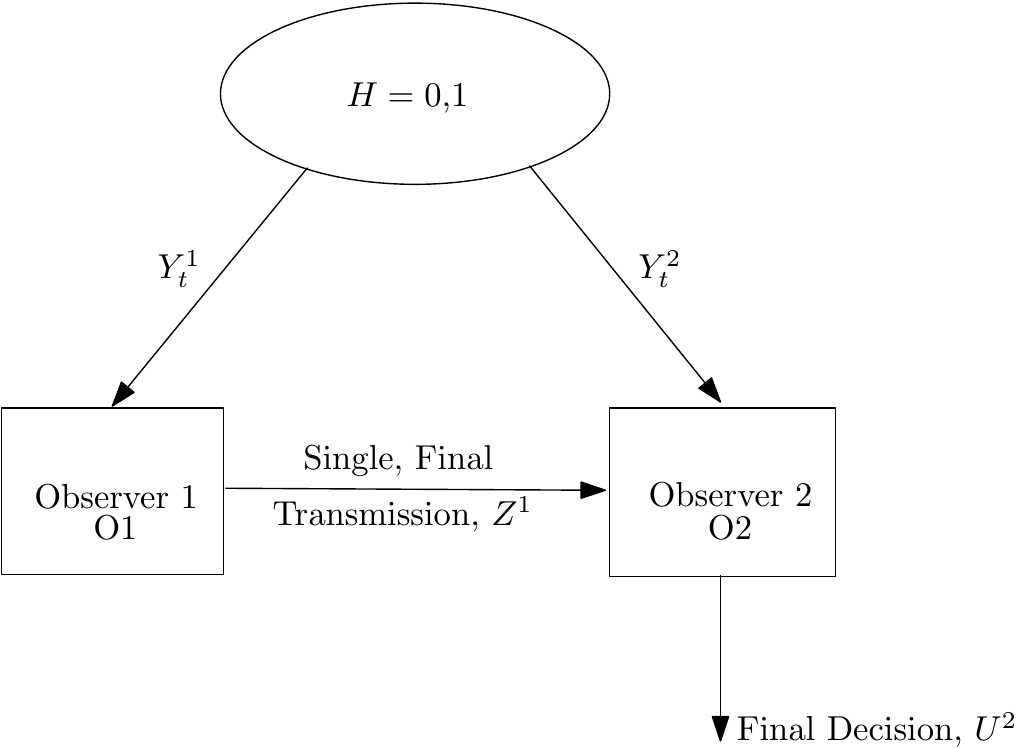}
\caption{Decentralized Detection}

\end{center} 
\end{figure}
Observer 1 observes the measurement process \(Y^1_t\), \(t=1,2,...\). At any time \(t\), after having observed the sequence of observations \(Y^1_{1:t}\), observer 1 can decide either to stop and send a binary message 0 or 1 to observer 2 or to postpone its decision and get another measurement. Each time observer 1 postpones its decision, a cost of \(c^1\) is incurred. (The cost \(c^1\) incorporates the additional cost of taking a new measurement, the energy cost of staying on for another time step and/or a penalty for delaying the decision.) Note that observer 1 transmits only a single final binary message to observer 2. The decision of observer 1 at time \(t\) is based on the entire sequence of observations till that time, in other words, observer 1 has perfect recall. Thus, we have that 
\begin{equation}
 Z^1_t=\gamma^1_t(Y^1_{1:t}),
\end{equation}
where \(Z^1_t\) is observer 1's  message at time \(t\) to observer 2 and \(\gamma^1_t\) is the \emph{decision-function} used by O1 at time \(t\). \(Z^1_t\) belongs to the set \(\{0,1,b\}\), where we use \(b\) for blank message, that is, no transmission. The sequence of functions \(\gamma^1_t,t=1,2,...,\) constitute the \emph{policy} of observer 1. Let \(\tau^1\) be the stopping time when observer sends a final message to observer 2, that is,
\begin{equation} \label{eq:Pf1}
 \tau^1=min\{t: Z^1_t \in \{0,1\}\}
\end{equation}

We allow two possibilities for the operation of observer 2. \\
\emph{Case A}: In this case, O2 first waits for O1 to send a final message. After receiving observer 1's final message, observer 2 can decide either to stop and declare a decision on the hypothesis or take additional measurements on its own.  After observer 2 has made \(k\) measurements (\(k=1,2,...\)), it can decide to stop and declare a final decision on the hypothesis or take a new measurement. Each time observer 2 decides to take another measurement it incurs a cost \(c^2\). Whenever observer 2 makes a final decision \(U \in \{0,1\}\) on the hypothesis, it incurs a cost \(J(U,H)\). As in the case of observer 1, we assume observer 2 has perfect recall. Let \(U^2_k \in \{0,1,N\}\) be the decision made by observer 2 after receiving \(\tau^1\) messages \((Z^1_{1:\tau^1})\) from observer 1 and subsequently making \(k\) observations \((Y^2_{1:k})\) of its own, (where we use \(N\) for a null decision, that is, a decision to continue taking measurements). Thus,
\begin{equation} \label{eq:Pf2}
 U^2_k=\gamma^2_k(Y^2_{1:k},Z^1_{1:\tau^1}),
\end{equation}
for \(k=0,1,2,\ldots\). The message sequence  \(Z^1_{1:\tau^1}\) is a sequence of \(\tau^1-1\) blank messages followed by \(Z^1_{\tau^1} = 0\) or \(1\). The sequence of decision-functions \(\gamma^2_k,k=0,1,2,...\) constitute the \emph{policy} of observer 2. We define \(\tau^2\) to be the number of  measurements taken before observer 2 announces its final decision on the hypothesis, that is,
\begin{equation} \label{eq:Pf3}
 \tau^2=min\{k: U^2_k \in \{0,1\}\}
\end{equation}
\par
\emph{Case B}: In this case, O2 starts taking measurements at time \(t=1\) without waiting for O1 to send a final message. At time $t=1,2,\ldots,$ we have the following time-ordering of the two observers' observations and decisions:
\begin{figure}[H]
\begin{center}

\includegraphics[height=2cm,width=6cm]{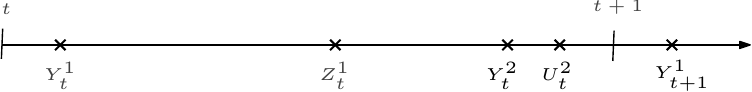}
\caption{Time ordering in P2}
\end{center}

\end{figure}
Thus, observer 2's decision at time \(t\) can be described as:
\begin{equation} \label{eq:Pf4}
 U^2_t=\gamma^2_t(Y^2_{1:t},Z^1_{1:t})
\end{equation}  
where \(U^2_t \in \{0,1,N\}\). This decision is a function of the observations made at O2 (\(Y^2_{1:t}\)) and the messages received from O1 (\(Z^1_{1:t}\)). The message sequence \(Z^1_{1:t}\) could be a sequence of \(t\) blank messages received from O1 or \(k\) blanks \((k<t)\) followed by a 0 or 1. Let \(\tau^2\) be the stopping time when observer 2 announces its final decision on the hypothesis, that is,
\begin{equation} \label{eq:Pf5}
 \tau^2=min\{t: U^2_t \in \{0,1\}\}
\end{equation}
 Note that we allow O2 to declare a final decision without getting the final message from O1. Also, O1 does not know whether O2 has stopped or not, that is, there is no feedback from O2 to O1. As in Case A, a penalty of $c^2$ is incurred every time O2 decides to postpone its final decision and a terminal cost of $J(U,H)$ is incurred when O2 makes its final decision $U \in \{0,1\}$.  
 \par 
 In both the cases above, we assume that the cost parameters \(c^1,c^2\) are finite positive numbers and \(J(U,H)\) is non-negative and bounded by a finite constant \(L\) for all \(U\) and \(H\). Moreover, we assume that cost of an error in the final decision is more than cost of a correct decision, that is, \( J(0,1) > J(1,1)\) and \(J(1,0) > J(0,0)\).
 We can now formulate an optimization problem for each of the two cases above.
 \subsubsection {Problem P1} 
 We consider a finite horizon \(T^1\) for observer 1. That is, if the observer 1 has not sent its final message till time \(t=T^1-1\), it must do so at time \(T^1\). In other words, we require that \(\tau^1 \leq T^1\).
 Similarly for observer 2 described in Case A above, we require that it can at most take 
 \(T^2\) measurements before declaring its final decision, that is, \(\tau^2 \leq T^2\). The optimization problem is to select policies \(\Gamma^1 = (\gamma^1_1,\gamma^1_2,...,\gamma^1_{T^1})\) and \(\Gamma^2 = (\gamma^2_0,\gamma^2_1,\gamma^2_2,...,\gamma^2_{T^2})\) to minimize
 \begin{equation} \label{eq:Objective1}
 \mathds{E}^{\Gamma^1,\Gamma^2}\{c^1\tau^1+c^2\tau^2+J(U^2_{\tau^2},H)\}
 \end{equation}
 where \(\tau^1,\tau^2\) and \(U^2_k\), \(k=0,1,\ldots\) are defined by equations (\ref{eq:Pf1}), (\ref{eq:Pf2}) and (\ref{eq:Pf3}) above. 
 
 \subsubsection{Problem P2}
   As in Problem P1, we have a finite horizon \(T^1\) for O1, that is, \(\tau^1 \leq T^1\) and a finite horizon \(T^2 (\geq T^1)\) for O2. O2's operation is as described in Case B above. The optimization problem is to select policies \(\Gamma^1 = (\gamma^1_1,\gamma^1_2,...,\gamma^1_{T^1})\) and \(\Gamma^2 = (\gamma^2_1,\gamma^2_2,...,\gamma^2_{T^2})\) to minimize
 \begin{equation}\label{eq:Objective2}
 \mathds{E}^{\Gamma^1,\Gamma^2}\{c^1\tau^1+c^2\tau^2+J(U^2_{\tau^2},H)\}
 \end{equation}
 where \(\tau^1,\tau^2\) and \(U^2_t\), \(t=1,2,\ldots\) are defined by equations (\ref{eq:Pf1}), (\ref{eq:Pf4}) and (\ref{eq:Pf5}) above. 
 
 \subsection{Features of the Problem}
 In both the problems formulated above, the two observers share a common system objective given by equations (\ref{eq:Objective1}) or (\ref{eq:Objective2}). The two observers, however, make decisions based on different information. Thus, Problems P1 and P2 are \emph{team problems}. Moreover, since the actions of observer 1 influence the information available to observer 2, these are dynamic team problems \cite{Ho1980}.
Dynamic team problems are known to be hard as they usually involve non-convex functional optimization over the space of policies of the decision-makers. Finding structural results for these problems is an important step toward reducing the complexity of these problems. In the next two sections, we will establish qualitative properties of the optimal policies of the two observers. 

\section{Qualitative Properties of Optimal Policies for Observer O1} \label{sec:QP1}

\subsection{Information state for O1}
     Consider Problem P2 first. We first derive an information state for O1. For that purpose, we define:
     \begin{equation}
      \pi^1_t(Y^1_{1:t}) := P(H=0|Y^1_{1:t})
     \end{equation}
    The probability \(\pi^1_t\) is observer 1's belief on the hypothesis based on its sequence of observations till time \(t\). (For \(t=0\), we have \(\pi^1_0 = p_0\)). The following result provides a characterization of O1's optimal policy.
\begin{theorem}
  For Problem P2, with an arbitrary but fixed policy \(\Gamma^2\) of O2, there is an optimal policy for O1 of the form:
  \begin{equation} \label{eq:Qp1}
  Z^1_t = \gamma^1_{t}(\pi^{1}_{t})
  \end{equation}
  for \(t=1,2,...,T^1\).
  In particular, if globally optimal policies \(\Gamma^{1,*},\Gamma^{2,*}\) exist, then \(\Gamma^{1,*}\) can be assumed to be of the from in (\ref{eq:Qp1}) without loss of optimality. Moreover, for a fixed \(\Gamma^2\), the optimal policy of O1 can be determined by selecting the minimizing option at each step of the following dynamic program:
  \begin{align}
    &V_{T^1}(\pi) := min \{ \nonumber \\ &\mathds{E}^{\Gamma^2}[c^2\tau^2+J(U^2_{\tau^2},H)|\pi^1_{T^1}=\pi,Z^1_{1:T^1-1}=b_{1:T^1-1},\nonumber \\&\hspace{20pt}Z^1_{T^1}=0],  \nonumber \\
                         &\mathds{E}^{\Gamma^2}[c^2\tau^2+J(U^2_{\tau^2},H)|\pi^1_{T^1}=\pi,Z^1_{1:T^1-1}=b_{1:T^1-1},\nonumber \\&\hspace{20pt}Z^1_{T^1}=1]  \}
 \end{align}   
 and for \(k=(T^1-1),...,2,1\),
 \begin{align}
  &V_{k}(\pi) := min\{ \nonumber \\ &\mathds{E}^{\Gamma^2}[c^2\tau^2+J(U^2_{\tau^2},H)|\pi^1_{k}=\pi,Z^1_{1:k-1}=b_{1:k-1},Z^1_{k}=0], \nonumber \\
                          &\mathds{E}^{\Gamma^2}[c^2\tau^2+J(U^2_{\tau^2},H)|\pi^1_{k}=\pi,Z^1_{1:k-1}=b_{1:k-1},Z^1_{k}=1], \nonumber \\
                          &c^1 + \mathds{E}^{\Gamma^2}[V_{k+1}(\pi^1_{k+1})|\pi^1_{k}=\pi,Z^1_{1:k}=b_{1:k}]                               \}
 \end{align}
 where the superscript \(\Gamma^2\) in the expectation denotes that the expectation is defined for a fixed choice of \(\Gamma^2\). (\(Z^1_{1:k}=b_{1:k}\) denotes a sequence of \(k\) blank messages.)
\end{theorem} 
\begin{proof}
  See Appendix A.
  \\
\end{proof}
\par
   The result of Theorem 1 can be intuitively explained as follows. At any time \(t\), if the observer 1 has not already sent its final message, it has to choose between three choices of action - send \(0,1\) or \(b\). In order to evaluate the expected cost of sending a \(0\) or \(1\), O1 needs a belief on the state of the environment, that is, a belief on \(H\) and a belief on the information available to O2. Since O1 has not yet sent a final message, the information at O2 consists of \(Z^1_{1:t-1} = b_{1:t-1}\), the decision of O1 at time \(t\) \((Z^1_t)\) and the observations that O2 has made or may make in the future. Thus O1 needs to form a belief on \(Y^2_{1:T^2}\), since the rest of O2's information it already knows. Now because of conditional independence of observations at the two observers, it is sufficient to form a belief on \(H\)  to know the probabilities of \(Y^2_{1:T^2}\). Similarly, to evaluate the cost of sending a \(b\), O1 needs to form a belief on O2's information and what information O1 may obtain by future measurements - \(Y^1_{t+1:T^1}\). Once again, conditional independence of the observations made at different times given \(H\) implies a belief on \(H\) is sufficient to evaluate the cost of this action as well. These arguments indicate that the decisions at O1 should be made based only on its belief on \(H\), that is, \(\pi^1_t\).
\par
\emph{Corollary:} Theorem 1 holds for O1 in Problem P1 also.
\begin{proof}
This result can be obtained by following the steps in Appendix A without any modifications.
An intuitive explanation of this result is as follows: In the proof of Theorem 1, we fixed \(\Gamma^2\) to any arbitrary choice. In particular, consider any policy of O2 that waits till it gets a final decision from O1. After it receives the final message from O1 at time \(\tau^1\), it uses only observations made after \(\tau^1\) to make a decision. This class of policies is essentially the policies available to O2 in problem P1. Since the optimal structure of O1's policy as given in (\ref{eq:Qp1}) holds for any choice of \(\Gamma^2\), it also holds for all possible policies of O2 in problem P1. 
\end{proof} 


\subsection{Classical Two-Threshold Rules Are Not Optimal}
      In the sequential detection problem with a single observer \cite{Wald}, it is well known that an optimal policy is a function of the observer's belief \(\pi_t\) and is described by two thresholds at each time \(t\). That is the decision at time \(t\), \(Z_t\) is given as:
      \[ Z_t = \left \{ \begin{array}{ll}
               1 & \mbox{if $\pi_t \leq \alpha_t$} \\
               N & \mbox{if $\alpha_t<\pi_t < \beta_t$} \\
               0 & \mbox{if $\pi_t \geq \beta_t$}
               \end{array}
               \right. \]
     where \(N\) denotes a decision to continue taking measurement and \(\alpha_t \leq \beta_t\) are real numbers in \([0,1]\).
 A similar two-threshold structure of optimal policies was also established for the decentralized Wald problem in \cite{Dec_Wald}. We will show by means of two counterexamples that such a structure is not necessarily optimal for observer 1 in Problem P1. Example 2 is similar to an example demonstrating the sub-optimality of threshold rules in a more general decentralized sequential detection problem that appeared in \cite{Tsitsiklis_86}. 
 \par
 \emph{Example 1} 
 \par
 Consider the following instance of Problem P1. We have equal prior on \(H\), that is \(P(H=0)=P(H=1)=1/2\). O1 has a time horizon of \(T^1=2\). Its observation space is \(\mathcal{Y}^1 = \{1,2,3\}\). The observations at time \(t=1\) have the following conditional probabilities:
 \[ \begin{array}{lclc}
Observation & 1 & 2 & 3 \\
P(.|H=0)    & 0 & p & (1-p) \\
P(.|H=1)    & (1-p) & p & 0 
\end{array}\]
and at time \(t=2\) have the following conditional probabilities:
\[ \begin{array}{lclc}
Observation & 1 & 2 & 3 \\
P(.|H=0)    & 0 & q & (1-q) \\
P(.|H=1)    & (1-q) & q & 0 
\end{array}\]
where \(p,q \in [0,1]\).
Observe that O1's belief on \(\{H=0\}\) (that is, \(\pi^1\)), only takes 3 possible values - \(0,1\) and \(1/2\) after any number of measurements.
O1 has to send a final message - \(0\) or \(1\) - to O2 by time \(T^1=2\). If O1 delays sending its final message to time \(t=2\), an additional cost \(c^1\) is incurred. After receiving a message from O1, observer 2 can either declare a decision on the hypothesis or take at most 1 more measurement of its own, that is, we have \(T^2 =1\). The measurements of O2 are assumed to be noiseless, so when O2 takes a measurement it knows exactly the value of \(H\). However, the measurement comes at a cost of \(c^2\). We assume that \(J(U,H) =0\) if \(U=H\), and in the case of  a mistake \((U \neq H)\), we assume that the cost is sufficiently high so that unless O2 is certain from O1's messages what the true hypothesis is, it will prefer taking a measurement at a cost \(c^2\) rather than making a guess. 
  At \(p=0.6\), \(c^2 > 3c^1\), it can be easily verified that the best threshold rule for observer 1 is described as follows:
  \[ Z^1_1 = \left \{ \begin{array}{ll}
               1 & \mbox{if $\pi^1_1 = 0$} \\
               b & \mbox{if $\pi^1_1 = 1/2$} \\
               0 & \mbox{if $\pi^1_1 = 1$}
               \end{array}
               \right. \]
               and
               \[ Z^1_2 = \left \{ \begin{array}{ll}
               1 & \mbox{if $\pi^1_1 = 0$} \\
               0 & \mbox{if $\pi^1_1 > 0$}
               \end{array}
               \right. \]
  
  If observer 2 receives \(0\) or \(1\) at time \(t=1\), it declares the received message as the final decision on the hypothesis, otherwise it waits for a final message from O1. At \(t=2\), if O2 receives \(1\), it declares \(1\) as the final decision, otherwise it takes a measurement. Then the expected cost of this policy is given as: \( pc^1 + p(1+q)c^2/2\) (since the system incurs a cost \(c^1\) with probability \(p\) and a cost \(c^2\) with probability \(p/2 + pq/2\)).
  \par
  Now consider the following non-threshold policy for observer 1,
  \[ Z^1_1 = \left \{ \begin{array}{ll}
               1 & \mbox{if $\pi^1_1 = 0$} \\
               0 & \mbox{if $\pi^1_1 = 1/2$} \\
               b & \mbox{if $\pi^1_1 = 1$}
               \end{array}
               \right. \] 
               and
               \[ Z^1_2 = \left \{ \begin{array}{ll}
               1 & \mbox{if $\pi^1_1 = 0$} \\
               0 & \mbox{if $\pi^1_1 > 0$}
               \end{array}
               \right. \]
  
  Unlike a classical two-threshold rule, the above rule requires O1 to send a blank symbol at time $1$ even though O1 is certain that true $H$ is $0$. If observer 2 receives \(0\) at time \(t=1\), it takes a measurement and incurs a cost \(c^2\). If O2 receives a \(1\) at \(t=1\), it declares \(1\) as the final decision. If O2 receives a \(b\) at time \(t=1\), it waits for the final message at \(t=2\) and then declares the received message as its final decision on the hypothesis.  Then the expected cost of this policy is given as: \( pc^2 + (1-p)c^1/2\) (since the system incurs a cost \(c^2\) with probability \(p\) and a cost \(c^1\) with probability \((1-p)/2\)). It is now easily seen that at \(p=0.6\) and \(c^2>3c^1\), if we choose \(q > 1-\frac{4c^1}{3c^2}\), the non-threshold policy outperforms the best threshold policy.

 \emph{Discussion of the Example:} The principle behind a threshold rule is to stop and send a message if O1 is certain, otherwise postpone the decision and take another measurement. The additional cost of delay is justified by the likelihood of getting a good measurement in the next time instant. In our example, if O1 gets the observations \(1\) or \(3\) at \(t=1\) and is able to convey to O2 that it is certain about the true hypothesis and what this hypothesis is, then it prevents O2 from taking a measurement thus saving a cost \(c^2\). The threshold rule achieves this objective by sending \(0\) for observation \(3\) and \(1\) for observation \(1\). However, in the case when O1 gets measurement \(2\), it decides to wait for the next observation. By choosing \(q\) sufficiently high, the likelihood of getting a good measurement at \(t=2\) can be made very low. In this case, the cost of delaying a decision \((c^1)\) begins outweighing the expected payoff from a new measurement. The non-threshold rule essentially tries to correct this drawback. If at time \(t=1\), O1 gets measurement \(2\), it stops and sends \(0\) to O2. At O2, this is interpreted as a message to go and take measurement of its own. Note that the non-threshold rule still ensures that whenever O1 is certain about \(H\), it is able to send enough information to O2 to prevent it from taking a measurement.  
 \par
 \emph{Example 2}
 \par
  Consider the same problem as in Example 1 but with O1's observations at \(t=1\) now given by the following conditional probabilities.  
\[ \begin{array}{lcccc}
Observation & 1 & 2 & 3 & 4 \\
P(.|H=0)    & 0 & p/3 & 2p/3 & (1-p) \\
P(.|H=1)    & (1-p) & 2p/3 & p/3 & 0 
\end{array}\]
O1's observations at time \(t=2\) are just noise and give no new information. The rest of the model is same as in Example 1. Note that the observations are indexed in order of the posterior belief \(\pi^1\) they generate, that is, \(P(H=0|Observation 1) < P(H=0|Observation 2)\) and so on. If O1 postpones its final message to time \(t=2\), it has to pay an additional cost of \(c^1\). Observer 2 can make a noiseless measurement at a cost of \(c^2\). As in Example 1, Observer 2's cost of making a wrong decision is chosen sufficiently high so that unless it is certain from O1's message what the true hypothesis is, O2 will prefer taking a measurement at a cost \(c^2\) than making a guess. 
It can be shown that for equal prior (\( p_0 =1/2\)), \(c^2>2c^1\) and \(1/2<p<1\), a non-threshold rule for O1 (given below) performs better than any threshold policy. 
\begin{itemize}
\item At \(t=1\), send \(0\) if observation 2 occurs and \(1\) if observation 3 occurs. Send a blank otherwise.
\item At \(t=2\), send \(1\) if \(\pi^1_2\) is less than \(1/2\) and \(0\) otherwise.
\end{itemize} 

The corresponding policy for O2 is as follows:
\begin{itemize}
\item At \(t=1\), if a \(0\) or \(1\) is received, take a measurement, otherwise wait till $t=2$.
\item At \(t=2\), declare the receive symbol as the final decision.
\end{itemize}
The cost of the above choice of policies is: $pc^2 + (1-p)c^1$.
 \subsection{Parametric Characterization of Optimal Policies}
   An important advantage of the threshold rules in the case of the centralized or the decentralized Wald problem is that it modifies the problem of finding the globally optimal policies from a sequential functional optimization problem to a sequential parametric optimization problem. Even though we have established that a classical two-threshold rule does not hold for our problem, it is still possible to get a finite parametric characterization of an optimal policy for observer 1. Such a parametric characterization provides significant computational advantage in finding optimal policies, for example by reducing the search space for an optimal policy.         \par
 In Theorem 1, we have established  that for an arbitrarily fixed choice of O2's policy, the optimal policy for O1 can be determined by backward induction using the functions \(V_k(\pi), k=T^1,...,2,1\). We will call \(V_k\) the value function at time \(k\).
 We have the following lemma.
 \begin{lemma}
  In problem P1 or P2, with a fixed (but arbitrary) choice of \(\Gamma^2\), the value function at \(T^1\) can be expressed as:
  \begin{equation}
  V_{T^1}(\pi) := min\{ L^{0}_{T^1}(\pi), L^{1}_{T^1}(\pi)\}
  \end{equation}
  where \(L^{0}_{T^1}(\cdot)\) and \(L^{1}_{T^1}(\cdot)\) are affine functions of \(\pi\) that depend on the choice of O2's policy \(\Gamma^2\).
  Also, the value function at time \(k\) can be expressed as:
  \begin{equation}\label{eq:PFlemma1}
  V_{k}(\pi) := min\{ L^{0}_k(\pi), L^{1}_k(\pi), G_k(\pi)\}
  \end{equation}
  where \(L^{0}_{k}(\cdot)\) and \(L^{1}_{k}(\cdot)\) are affine functions of \(\pi\) and \(G_k(\cdot)\) is a concave function of \(\pi\). The functions \(L^{0}_{k}(\cdot)\), \(L^{1}_{k}(\cdot)\) and \(G_k(\cdot)\) depend on the choice of \(\Gamma^2\).
 \end{lemma}
 \begin{proof}
 See Appendix B.
\end{proof}
 
 \begin{theorem}
    In Problem P1 or P2, for any fixed policy \(\Gamma^2\) of O2, an optimal policy for O1 can be characterized by \emph{at most 4 thresholds}. In particular, without any loss in performance, one can assume O1's policy to be of the following form:
   \[ Z^1_{T^1} = \left \{ \begin{array}{ll}
               1 & \mbox{if $\pi^1_{T^1} \leq \alpha_{T^1}$} \\
               0 & \mbox{if $\pi^1_{T^1} > \alpha_{T^1}$ }
               \end{array}
               \right. \]
               where \(0 \leq \alpha_{T^1} \leq 1\)  and for \(k=1,2,..,T^1-1\),
         \[ Z^1_k = \left \{ \begin{array}{ll}
               b & \mbox{if $\pi^1_k < \alpha_k$}  \\
               1 & \mbox{if $\alpha_k \leq \pi^1_k \leq \beta_k$} \\
               b & \mbox{if $\beta_k<\pi^1_k < \delta_k$} \\
               0 & \mbox{if $\delta_k \leq \pi^1_k \leq \theta_k$} \\
               b & \mbox{if $\pi^1_k > \theta_k$}
               \end{array}
               \right. \]  
               where \(0 \leq \alpha_k \leq \beta_k \leq \delta_k \leq \theta_k \leq 1\). 
\end{theorem}
\begin{proof}
Theorem 2 is an immediate consequence of Lemma 1, since taking minimum of straight lines and concave functions can partition the interval \( [0,1]\) into at most  five regions. The thresholds above essentially signify the boundaries of these regions. For a given \(\Gamma^2\), it is possible that at some time instant \(k\), the optimal policy for O1 partitions the belief interval \([0,1]\) as \(\{b,0,b,1,b\}\) instead of \(\{b,1,b,0,b\}\). In this case, it is easily seen that simply interchanging the roles of \(0\) and \(1\) in O1's policy and in \(\Gamma^2\) at time \(k\) would result in the threshold structure of the theorem without loss of performance.
\end{proof}
 It is of course possible that in specific cases, some of these five regions are absent which would correspond to some of the above thresholds having the same value. For example, in the non-threshold rule given in the Example 1 earlier, the rule at \(t=1\) corresponds to having \(\alpha=0\) and 
\(\beta=\delta\) which results in a 3-interval partition of \([0,1]\) corresponding the rule given there.

\section{Qualitative Properties for Observer O2} \label{sec:QP2}

\subsection{Problem P1}
Consider a fixed policy \(\Gamma^1 = (\gamma^1_1,\gamma^1_2,...,\gamma^1_{T^1})\) for O1. Then, after O1 sends its final message, we can define the following probability for O2:
\[\pi^2_0 := P^{\Gamma^1}(H=0|Z^1_{1:\tau^1}) \]
This is O2's belief on the true hypothesis after having observed the messages from O1 (that is a sequence of \(\tau^1-1\) blanks and a final \(Z^1_{\tau^1} \in \{0,1\}\)). Now, the optimization problem for O2 is the classical centralized Wald problem  \cite{Wald} with the prior probability given by \(\pi^2_0\). It is well-known that the optimal policy for the Wald problem is a rule of the form:
\[\hspace{20pt} U^2_k = \left \{ \begin{array}{ll}
               1 & \mbox{if $\pi^2_k \leq w^1_k$} \\
               N & \mbox{if $w^1_k<\pi^2_k<w^2_k$} \\
               0 & \mbox{if $\pi^2_k \geq w^2_k$}
               \end{array}
               \right. \]               
   where \(\pi^2_k\) is the belief on hypothesis after \(k\) observations, 
   \[ \pi^2_k(Y^2_{1:k}) := P^{\Gamma^1}(H=0|Z^1_{1:\tau^1},Y^2_{1:k})\] 
   \[ = \frac{P(Y^2_{1:k}|H=0).\pi^2_0}{P(Y^2_{1:k}|H=0).\pi^2_0+P(Y^2_{1:k}|H=1).(1-\pi^2_0)}, \]
and \(w^1_k \leq w^2_k\), for \(k=0,1,2,..,T^2-1\) and \(w^1_{T^2} = w^2_{T^2}\) are the optimal thresholds for the Wald problem with horizon \(T^2\).   

\subsection{Information State in Problem P2}
  Consider a fixed policy \(\Gamma^1 = (\gamma^1_1,\gamma^1_2,...,\gamma^1_{T^1})\) for O1. Define the following probability for O2:
  \[ \pi^2_t(Y^2_{1:t},Z^1_{1:t}) := P^{\Gamma^1}(H=0|Y^2_{1:t},Z^1_{1:t}) \]
   \(\pi^2_t\) is observer 2's belief on the hypothesis based on its observations till time \(t\) and the messages received from O1 till time \(t\) (where the messages from O1 could be all blanks or some blanks terminated by a \(0\) or \(1\)). For \(t=0\), we have \(\pi^2_0 = p_0\).\\ 
The following theorem shows that \(\pi^2_t\) and \(Z^1_{1:t}\) together form an information state for O2.
 \begin{theorem}
  In Problem P2, with an arbitrary but fixed policy \(\Gamma^1\) of O1, there is an optimal policy for O2 of the form:
  \begin{equation} \label{eq:Qp21}
  U^2_t = \gamma^2_t(Z^1_{1:t},\pi^2_t)
  \end{equation}
  for \(t=1,2,...,T^2\). Moreover, this optimal policy can be determined by the following dynamic program:
  \begin{align}
  	\tilde{V}_{T^2}(z^1_{1:T^1},\pi) := min \{ &E^{\Gamma^1} [J(0,H)|\pi^2_{T^2}=\pi], \nonumber \\
  	                                  &E^{\Gamma^1} [J(1,H)|\pi^2_{T^2}=\pi]  \} \label{eq:Qp22.1}
  \end{align}
  and for \(k=(T^1-1),...,1\),
  \begin{align}
   &\tilde{V}_k(z^1_{1:k},\pi) := min \{ \nonumber \\
   &E^{\Gamma^1} [J(0,H)|\pi^2_{k}=\pi], \nonumber \\
  	                                  &E^{\Gamma^1} [J(1,H)|\pi^2_{k}=\pi], \nonumber \\
  	                                  & c^2 + E^{\Gamma^1} [\tilde{V}_{k+1}(Z^1_{1:k+1},\pi^2_{k+1})|\pi^2_{k}=\pi,Z^1_{1:k}=z^1_{1:k}]\} \label{eq:Qp22.2}
  \end{align}
 \end{theorem}    
\begin{proof}
  See Appendix C.
\end{proof}
\par
 Observe that in the last term of (\ref{eq:Qp22.2}), which corresponds to the cost of postponing the final decision at time \(k\), we have \(\pi^2_k\) as well as all messages from O1 in the conditioning variables. It is because of this term  that we need the entire sequence of messages as a part of the information state. To intuitively see why these messages are needed in the conditioning, note that the cost of continuing depends on future messages from O1. In order to form a belief on future messages, O2 needs a belief on the hypothesis and (since  O1 has perfect recall) a belief on all observations of O1 so far.  Clearly, the messages received till time \(k\) provide information about the observations of O1 till time \(k\) and are therefore included in the information state.
 \par
 We can now prove the following lemma about the value functions \(\tilde{V}_k\).
 \begin{lemma}
 The value function at \(T^2\) can be expressed as:
  \begin{equation}
  \tilde{V}_{T^2}(z^1_{1:T^1},\pi) := min\{ l^{0}(\pi), l^{1}(\pi)\}
  \end{equation}
  where \(l^{0}\) and \(l^{1}\) are affine functions of \(\pi\) that are independent of the choice of O1's policy \(\Gamma^1\).
  Also, the value function at time \(k\) can be expressed as:
  \begin{equation}
  \tilde{V}_{k}(z^1_{1:k},\pi) := min\{ l^{0}(\pi), l^{1}(\pi), G_k(z^1_{1:k},\pi)\}
  \end{equation}
   where, for each realization \(z^1_{1:k}\) of messages from O1, \(G_k\) is a concave function of \(\pi\) that depends on the choice of O1's policy, \(\Gamma^1\).
 \end{lemma}
 \begin{proof}
 See Appendix D.
 \end{proof}
 \begin{theorem}
   For a fixed policy \(\Gamma^1\) of O1, an optimal policy of O2 is of the form:
    \[ U^2_{T^2} = \left \{ \begin{array}{ll}
               1 & \mbox{if $\pi^2_{T^2} \leq \alpha_{T^2}$} \\
               0 & \mbox{if $\pi^2_{T^2} > \alpha_{T^2}$ }
               \end{array}
               \right. \]
   \[ U^2_k = \left \{ \begin{array}{ll}
               1 & \mbox{if $\pi^2_k \leq \alpha_k(Z^1_{1:k})$}  \\
               N & \mbox{if $\alpha_k(Z^1_{1:k}) < \pi^1_k < \beta_k(Z^1_{1:k})$} \\
               0 & \mbox{if $\pi^1_k \geq \beta_k(Z^1_{1:k})$}
               \end{array}
               \right. \] 
               where \(0 \leq \alpha_k(Z^1_{1:k}) \leq \beta_k(Z^1_{1:k}) \leq 1\) are thresholds that depend on sequence of messages received from O1 (\(Z^1_{1:k}\)).
 \end{theorem}
 \begin{proof}
     At any time \(k\), if \(\pi^2_k = 0\), then it is optimal to stop and declare the hypothesis to be \(1\) since cost of continuing will be at least \(c^2+J(1,1)\) which is more than \(J(1,1)\) - the cost of immediately stopping and declaring \(U^2_k =1\). Similarly, at \(\pi^2_k = 1\), it is optimal to stop and declare \(U^2_k=0\). These observations along with the fact that the value functions $\tilde{V}_k$ are minimum of affine and concave functions for each realization of the messages received imply the result of the theorem. 
 \end{proof}
          
          		Thus, according to Theorem 4, the thresholds to be used at time \(k\) by O2 depend on the sequence of messages received from O1 until time \(k\). This kind of parametric characterization may not appear very appealing since for each time \(k\) one may have to know a number of possible thresholds - one for each possible realization of messages \(z^1_{1:k}\). We will now argue that there is in fact a simple representation of the thresholds. Note that after time \(\tau^1\), when O1 sends a final message, O2 is faced with a classical Wald problem with an available time-horizon of \(T^2-\tau^1\). Now suppose that the classical Wald thresholds are available for a time horizon of length \(T^2\) -lets call these \((w^1_0,w^2_0),(w^1_1,w^2_1),(w^1_2,w^2_2),...,w_{T^2}\). Then the Wald thresholds for a problem with time horizon \(T^2-\tau^1\) are simply \((w^1_{\tau^1},w^2_{\tau^1}),(w^1_{\tau^1+1},w^2_{\tau^1+1}),(w^1_{\tau^1+2},w^2_{\tau^1+2}),...,w_{T^2}\). Thus, once  O2 hears a final message from O1, it starts using the classical Wald thresholds from that time onwards. In other words, O2 operation is described by the following simple algorithm: 
   \begin{itemize}
   \item From time \(k=1\) onwards, the optimal policy is to use a threshold rule given by 2 numbers \(\alpha_k(b_{1:k})\) and \(\beta_k(b_{1:k})\), \emph{until O1 sends its final message \(Z^1_k \in \{0,1\}\).} (As before, \(b_{1:k}\) stands for sequence of \(k\) blank messages.)
   \item If O1 sends the final message at time $k$, start using Wald thresholds: \((w^1_k,w^2_k),...,w_{T^2}\).
   \end{itemize}
          		   Thus O2's optimal policy is completely characterized by just two tables of thresholds: \([(\alpha_1(b_{1:1}),\beta_1(b_{1:1}))\), \((\alpha_2(b_{1:2}),\beta_2(b_{1:2})),...,(\alpha_{T^1}(b_{1:T^1}),\beta_{T^1}(b_{1:T^1}))]\) and the Wald thresholds \([(w^1_0,w^2_0),(w^1_1,w^2_1),(w^1_2,w^2_2),...,w_{T^2}]\) .

\section{Optimal Policies} \label{sec:Optimal_Policies}
 In the previous sections, we identified qualitative properties of the optimal policies for the two observers. Moreover, if the policy $\Gamma^2$ $(\Gamma^1)$ of O2 (O1) has been chosen already, Theorems 1 (Theorem 3) provides a dynamic programming solution to find an optimal policy $\tilde{\Gamma}^{1}$ of O1  ($(\tilde{\Gamma}^{2})$ of O2) for the given choice of $\Gamma^2$ $(\Gamma^1)$. An iterative application of such an approach may be used to identify person-by-person optimal pair of strategies. However, finding globally optimal or near optimal strategies for such dynamic team problems remains a challenging task since it involves non-convex functional optimization \cite{Ho1980}. In this section, we will  give a sequential decomposition of the global optimization problem. Such a decomposition provides a systematic methodology to find globally optimal or near-optimal policies for the two observers. 
\subsection{Sequential Decomposition for Problem P1}
 In Problem P1, observer 2 waits to receive a final message from observer 1 before it starts taking its measurements. After receiving the final message, observer 2 is faced with the centralized sequential detection problem studied by Wald. For the Wald problem, the thresholds characterizing the optimal policy and the cost of the optimal policy are known. For a Wald problem with horizon \(T\) and a starting belief \(\pi\) on the event \(\{H=0\}\), the cost of using the optimal Wald thresholds is a function of the belief \(\pi\) which we will denote by \(K^{T}(\pi)\). Since the Wald thresholds for observer 2 are known (or can be calculated as in \cite{Wald}), the designer's task in problem P1 is to find the best set of thresholds to be used by observer 1. Finding the best thresholds for all times \(t=1\) to \(T^1\) is a formidable optimization problem. Firstly, the system objective (equation (\ref{eq:Objective1})) is a complicated function of the thresholds selected for observer 1. Moreover, the objective must be optimized over the space of sequences of thresholds to be used from time \(t=1\) to \(T^1\). Below, we show that the optimization problem can in principle be solved in a sequential manner. In the resulting sequential decomposition, at each step the optimization is over the set of thresholds to be used at a single time instant instead of the space of sequences of thresholds from time 1 to $T^1$. Though the original optimization  problem is decomposed into several ``simpler'' optimization problems, each of these remain difficult nonetheless. We believe that the decomposed problems may be more amenable to approximation techniques.
 \par
 We first define the following: 
\begin{definition}
For \(t=1,2,...,T^1\) and a given choice of observer 1's decision functions from time instant \(1\) to \(t-1\), that is, (\( \Gamma^1_{t-1} = (\gamma^1_1,\gamma^1_2,...,\gamma^1_{t-1})\)), define

 \[ \xi_t := P^{\Gamma^1_{t-1}}(H,\pi^1_t|Z^1_{1:t-1}=b_{1:t-1}) \]
 For \(t=1,2,...,T^1\) and for a given choice of functions (\(\Gamma^1_t = (\gamma^1_1,\gamma^1_2,...,\gamma^1_{t})\)), define
\[ \eta_t[z^1_{t}] := P^{\Gamma^1_{t}}(H,\pi^1_t|Z^1_{1:t-1}=b_{1:t-1},Z^1_{t}=z^1_{t}) \]
where \(z^1_t \in \{0,1,b\}\).
\end{definition}
\begin{lemma} Consider any policy \(\gamma^1_t\), \(t=1,2,...,T^1\) for observer 1 that is characterized by 4 thresholds \((\alpha^1_t,\beta^1_t,\delta^1_t,\theta^1_t)\), for \(t=1,2,...,T^1-1\) and a threshold \(\alpha^1_{T^1}\) at time \(T^1\) (Theorem 2). Then,
\begin{enumerate}[i)]
\item There exist transformations \(Q^1_t\)  for \(t=1,2,...,T^1\) such that
\[\eta_t[z^1_t] = Q^1_t(\xi_t, \gamma^1_t,z^1_t)\]
for \(z^1_t \in \{0,1,b\}\), and 
\item There exist transformations \(Q^2_t\), \(t=1,2,...,T^1-1\) such that 
\[\xi_{t+1} = Q^2_t(\eta_t[b])\]
\end{enumerate}
\end{lemma}
\begin{proof}
We first prove the second part of the lemma.
By definition,
\begin{align}
&\xi_{t+1}(h,\pi^1) = P^{\Gamma^1_{t}}(H=h,\pi^1_{t+1}=\pi^1|Z^1_{1:t}=b_{1:t}) \nonumber \\
&= P^{\Gamma^1_{t}}(H=h,T_t(\pi^1_t,Y^1_{t+1})=\pi^1|Z^1_{1:t}=b_{1:t}) \label{eq:GoP1.1}
\end{align}
where we used the fact that O1's belief at time \(t+1\) is a function of its belief at time \(t\) and the observation at time \(t+1\), that is, \(\pi^1_{t+1} = T_t(\pi^1_t,Y^1_{t+1})\) (see Appendix A, equation (\ref{eq:belief_update_O1})). The right hand side of (\ref{eq:GoP1.1}) can further be written as:
\begin{align}
= \int_{y,\pi'}&[\mathbbm{1}_{T_t(\pi',y)=\pi^1} \nonumber \\ &.P^{\Gamma^1_{t}}(H=h,\pi^1_t=\pi',Y^1_{t+1}=y|Z^1_{1:t}=b_{1:t})] \nonumber \\
= \int_{y,\pi'}&[\mathbbm{1}_{T_t(\pi',y)=\pi^1}.P(Y^1_{t+1}=y|H=h) \nonumber \\ &.P^{\Gamma^1_{t}}(H=h,\pi^1_t=\pi'|Z^1_{1:t}=b_{1:t})] \nonumber \\
= \int_{y,\pi'}&[\mathbbm{1}_{T_t(\pi',y)=\pi^1}.P(Y^1_{t+1}=y|H=h).\eta_t[b](h,\pi')] \label{eq:expression_1}
\end{align}
The above integral is a function of \(\eta_t[b]\) and known observation statistics. Thus $\xi_{t+1} = Q^2_t(\eta_t[b])$, where $Q^2_t$ is given by the expression in (\ref{eq:expression_1}). \\
For the first part of the lemma, consider
\begin{align}
&\eta_t[b] = P^{\Gamma^1_{t}}(H=h,\pi^1_t=\pi^1|Z^1_{1:t}=b_{1:t}) \nonumber \\ 
&= \frac{P^{\Gamma^1_{t}}(H=h,\pi^1_t=\pi^1,Z^1_t=b|Z^1_{1:t-1}=b_{1:t-1})}{P^{\Gamma^1_{t}}(Z^1_t=b|Z^1_{1:t-1}=b_{1:t-1})} 
\end{align}
Under the 4-threshold rule for observer 1, \(Z^1_t =b\) if \(\pi^1_t \in \mathcal{C}_t\), where \(\mathcal{C}_t := [0,\alpha^1_t) \cup (\beta^1_t,\delta^1_t) \cup (\theta^1_t,1]\). Therefore, the above probability can be written as:
\begin{align}
&= \frac{\mathbbm{1}_{\pi^1 \in \mathcal{C}_t}.P^{\Gamma^1_{t}}(H=h,\pi^1_t=\pi^1|Z^1_{1:t-1}=b_{1:t-1})}{P^{\Gamma^1_{t}}(\pi^1_t \in \mathcal{C}_t|Z^1_{1:t-1}=b_{1:t-1})} \nonumber \\
&= \frac{\mathbbm{1}_{\pi^1 \in \mathcal{C}_t}.\xi_t(h,\pi^1)}{\int_{h,\pi'}\mathbbm{1}_{\pi' \in \mathcal{C}_t}\xi_t(h,\pi')}
\end{align}

The above equation is a function of \( \xi_t\) and the thresholds selected by \(\gamma^1_t\). Similar analysis holds for $\eta_t[0]$ and $\eta_t[1]$. This concludes the proof of the lemma.
\end{proof}

We can now present a sequential decomposition of problem P1.
\begin{theorem}
 For \(t=1,2,...,T^1-1\), there exist functions \(\mathcal{R}_{t}(\xi_{t},\alpha^1_{t},\beta^1_{t},\delta^1_t,\theta^1_t)\) and \(\mathcal{R}^*_t(\xi_{t})\) where

 \[  \mathcal{R}^{*}_{t}(\xi_t) = \inf_{\alpha^1_{t},\beta^1_{t},\delta^1_t,\theta^1_t}\mathcal{R}_{t}(\xi_{t},\alpha^1_{t},\beta^1_{t},\delta^1_t,\theta^1_t)\]
 and for \(t=T^1\), there exist functions \(\mathcal{R}_{T^1}(\xi_{T^1},\alpha^1_{T^1})\)and \(\mathcal{R}^*_{T^1}(\xi_{T^1})\) where 
 \[ \mathcal{R}^{*}_{T^1}(\xi_{T^1}) = \inf_{\alpha^1_{T^1}}\mathcal{R}_{T^1}(\xi_{T^1},\alpha^1_{T^1})\]
 such that the optimal thresholds can be evaluated from these functions as follows:
 \begin{enumerate}
 \item Note that \(\xi_1 := P(H,\pi^1_1)\) is fixed a priori and does not depend on any design choice. The optimal thresholds at \(t=1\) for O1 are given by optimizing parameters in the definition of \(\mathcal{R}^*_1(\xi_1)\). \\
 \item Once O1's thresholds at \(t=1\) are fixed, \(\eta_1[b]\) and hence \(\xi_2\) are fixed by lemma 3. The optimal thresholds for O1 at time \(t=2\) are given  by optimizing parameters in the definition of \(\mathcal{R}^*_2(\xi_2)\) \\
\item Continuing sequentially, \(\xi_t\) is fixed by the choice of past thresholds, and the optimal thresholds for O1 at time \(t\) are given  by optimizing parameters in the definition of \(\mathcal{R}^*_t(\xi_t)\). 
\end{enumerate}  
\end{theorem}
\begin{proof} 
 We will prove the result by backward induction.  
 Consider first the final horizon for O1: \(T^1\). Assume that a designer has already specified functions \(\gamma^1_1,\gamma^1_2,...,\gamma^1_{T^1-1}\) for O1.  The designer has to select a function to be used by O1 at time \(T^1\) in case the final message has not been already sent (that is, \(Z^1_{T^1-1}=b_{1:T^1-1}\)). By Theorem 2, this function is characterized by a single threshold \(\alpha^1_{T^1}\). For any choice of \(\alpha^1_{T^1}\), the future cost for the designer is \(K^{T^2}(\pi^2_0)\), where \(K^{T^2}(\cdot)\) is the cost of using optimal Wald thresholds with a time-horizon \(T^2\) and \(\pi^2_0\) is O2's belief on \(\{H=0\}\) after receiving \(Z^1_{1:T^1}\). The expected future cost for the designer can therefore be expressed as:
\begin{align}
&\mathds{E}\{c^2\tau^2+J(U_{\tau^2},H)|Z^1_{1:T^1-1}=b_{1:T^1-1}\} \nonumber\\
&=\mathds{E}\{K^{T^2}(\pi^2_0)|Z^1_{1:T^1-1}=b_{1:T^1-1}\} \nonumber\\
 &=K^{T^2}(P(H=0|Z^1_{T^1}=0,Z^1_{1:T^1-1}=b_{1:T^1-1})) \nonumber \\&\cdot P(Z^1_{T^1}=0|Z^1_{1:T^1-1}=b_{1:T^1-1}) \nonumber \\ 
&+K^{T^2}(P(H=0|Z^1_{T^1}=1,Z^1_{1:T^1-1}=b_{1:T^1-1}))\nonumber \\&\cdot P(Z^1_{T^1}=1|Z^1_{1:T^1-1}=b_{1:T^1-1})   \nonumber\\
\nonumber\\
&= K^{T^2}(P(H=0|Z^1_{T^1}=0,Z^1_{1:T^1-1}=b_{1:T^1-1}))\nonumber \\& \cdot P(\pi^1_{T^1} > \alpha^1_{T^1}|Z^1_{1:T^1-1}=b_{1:T^1-1})  \nonumber\\
&+K^{T^2}(P(H=0|Z^1_{T^1}=1,Z^1_{1:T^1-1}=b_{1:T^1-1}))\nonumber \\&\cdot P(\pi^1_{T^1} \leq \alpha^1_{T^1}|Z^1_{1:T^1-1}=b_{1:T^1-1}) \nonumber\\
\nonumber \\ 
 &=: L_{T^1}(\eta_{T^1}[0],\eta_{T^1}[1],\xi_{T^1},\alpha^1_{T^1}) \label{eq:GoP1.2}
\end{align}
where we used the fact that the probabilities in the arguments of \(K^{T^2}(\cdot)\) are marginals of \(\eta_{T^1}[0]\) \(\eta_{T^1}[1]\) respectively and the probabilities multiplying the functions \(K^{T^2}\) are marginals of \(\xi_{T^1}\). Using Lemma 3, we can write (\ref{eq:GoP1.2}) as
\begin{gather}
 \begin{aligned}
&L_{T^1}(Q^1_{T^1}(\xi^1_{T^1},0,\alpha^1_{T^1}),Q^1_{T^1}(\xi^1_{T^1},1,\alpha^1_{T^1}),\xi_{T^1},\alpha^1_{T^1}) \nonumber \\
&=: \mathcal{R}_{T^1}(\xi_{T^1},\alpha^1_{T^1})
 \end{aligned}
\end{gather}
Thus, for a fixed choice of functions \(\gamma^1_1,\gamma^1_2,...,\gamma^1_{T^1-1}\) used till time \(T^1-1\), the designer's future cost at $T^1$ ,if the final message was not sent before $T^1$, is a function of \(\xi_{T^1}\) (that is induced by the choice of the past decision functions) and the threshold  \(\alpha^1_{T^1}\) it selects at time \(T^1\). To find the best choice of threshold, the designer has to select \(\alpha^1_{T^1}\) to minimize \(\mathcal{R}_{T^1}(\xi_{T^1},\alpha^1_{T^1})\). Define
\[ \mathcal{R}^{*}_{T^1}(\xi_{T^1}) = \inf_{\alpha^1_{T^1}}R(\xi_{T^1},\alpha^1_{T^1})\]
For a given \( \xi_{T^1}\), the function \(\mathcal{R}^{*}_{T^1}\) describes the optimal future cost for the designer and the optimizing \(\alpha^1_{T^1}\) gives the best threshold. 
\par
Now assume that $\mathcal{R}^{*}_{t+1}(\xi_{t+1})$ describes the designer's optimal future cost from time $t+1$. At time $t$, if the past decision functions \(\gamma^1_1,\gamma^1_2,...,\gamma^1_{t-1}\) have been specified already, the designer's task is  to select thresholds \(\alpha^1_{t},\beta^1_{t},\delta^1_{t},\theta^1_{t}\) to be used by O1 at time \(t\). For a given choice of these thresholds, the future cost for the designer is \(K^{T^2}(\pi^2_0)\) if O1 sends a final message at \(t\). If a blank message is sent at \(t\), the designer will use the best threshold at the next time \(t+1\) and the future cost will be \(c^1 + \mathcal{R}^{*}_{t+1}(\xi_{t+1})\). The expected future cost for the designer is therefore given as:
\begin{align}
&\mathds{E}\{c^1(\tau^1-t)+ c^2\tau^2+J(U_{\tau^2},H)|Z^1_{1:t-1}=b_{1:t-1}\} \nonumber\\
&=K^{T^2}(P(H=0|Z^1_{t}=0,Z^1_{1:t-1}=b_{1:t-1})) \nonumber \\&\cdot P(Z^1_{t}=0|Z^1_{1:t-1}=b_{1:t-1}) \nonumber \\ 
&+K^{T^2}(P(H=0|Z^1_{t}=1,Z^1_{1:t-1}=b_{1:t-1}))\nonumber \\&\cdot P(Z^1_{t}=1|Z^1_{1:t-1}=b_{1:t-1}) \nonumber \\
 &+ [c^1+\mathcal{R}^{*}_{t+1}(\xi_{t+1})] \cdot P(Z^1_{t}=b|Z^1_{1:t-1}=b_{1:t-1}) \nonumber\\
 \\
&= K^{T^2}(P(H=0|Z^1_{t}=0,Z^1_{1:t-1}=b_{1:t-1})) \nonumber \\&\cdot P(\delta^1_{t}<\pi^1_{t}<\theta^1_{t}|Z^1_{1:t-1}=b_{1:t-1}) \nonumber \\
&K^{T^2}(P(H=0|Z^1_{t}=1,Z^1_{1:t-1}=b)) \nonumber \\&\cdot P(\alpha^1_{t}<\pi^1_{t}<\beta^1_{t}|Z^1_{1:t-1}=b_{1:t-1}) \nonumber \\
 &+[c^1 + \mathcal{R}^{*}_{t+1}(Q^2_{t}(\eta_{t}[b]))]  \cdot P(\pi^1_{t} \in \mathcal{C}_t|Z^1_{1:t-1}=b_{1:t-1}) \nonumber\\
 \\& =:  L_{t}(\eta_{t}[0],\eta_{t}[1],\eta_{t}[b],\xi_{t},\alpha^1_{t},\beta^1_{t},\delta^1_{t}, \theta^1_{t}) \label{eq:GoP1.3}
\end{align}
where we used the fact that the probabilities in the arguments of \(K^{T^2}(\cdot)\) are marginals of \(\eta_{t}[0]\) \(\eta_{t}[1]\) respectively and the probabilities multiplying the functions \(K^{T^2}\) and \(\mathcal{R}_{t+1}\) are marginals of \(\xi_{t}\). Using Lemma 3, we can write (\ref{eq:GoP1.3}) as a function of \(\xi_{t}\) (that is induced by the choice of past functions used till time \(t-1\)) and the thresholds selected at time \(t\):
 \begin{align}
  \mathcal{R}_{t}(\xi_{t},\alpha^1_{t},\beta^1_{t},\delta^1_{t},\theta^1_{t}) 
\end{align}
To find the best choice of threshold, the designer has to select \((\alpha^1_{t},\beta^1_{t},\delta^1_{t},\theta^1_{t})\) to minimize \(\mathcal{R}_{t}(\xi_{t},\alpha^1_{t},\beta^1_{t},\delta^1_{t},\theta^1_{t})\). Define
\begin{align}  \mathcal{R}^{*}_{t}(\xi_{t})  = \inf_{\substack{\alpha^1_{t},\beta^1_{t},\\\delta^1_{t},\theta^1_{t}}}\mathcal{R}_{t}(\xi_{t},\alpha^1_{t},\beta^1_{t},\delta^1_{t},\theta^1_{t})\end{align}
For a given \( \xi_{t}\), the function \(\mathcal{R}^{*}_{t}\) describes the optimal future cost for the designer and the optimizing thresholds are the best thresholds. 
The above analysis can be inductively repeated for all time instants. 
\par
The optimal thresholds can therefore be evaluated as follows: 
At \(t=1\), \(\xi_1\) is fixed  a priori, therefore one can use \(\mathcal{R}^*_1\) to find the best thresholds at time \(t=1\). Once these are selected, \(\xi_2\) can be found using Lemma 3 and one can use \(\mathcal{R}^*_2\) to find the best thresholds at time \(t=2\) and so on.

\end{proof} 

\par
\emph{Discussion:}
   The problem of choosing the optimal thresholds for observer 1 can be viewed as a sequential problem for the designer as follows: At each time $t$, the designer must specify the thresholds to be used by observer 1 in case the final message has not already been sent. In other words, at each time $t$, one can think that the designer is aware of the messages sent from O1 to O2 until $t-1$ and  in case these were only blanks, the designer must choose the thresholds to be used by O1 at time $t$. Thus, the designer is faced with a sequential optimization problem with a fixed temporal ordering of its decisions. Observe also that the designer has perfect recall: it knows all messages sent till time $t$. The designer, therefore, has a sequential problem with a \emph{classical information structure} \cite{Witsenhausen-separation}. The proof of Theorem 5 essentially describes the dynamic program for the designer's problem. The belief $\xi_t$ serves as the designer's information state and the  functions $R^*_t(\xi_t)$ are essentially the value functions of the dynamic program. This approach of introducing a designer with access to the common information between observers (that is, the information known to both observers: the messages from O1 to O2 in Problem P1) so as to convert a decentralized problem to one with classical information structure is illustrated and fully explained in \cite[Section IV]{NayyarTeneketzis:2008} for a communication problem. We refer the reader to that paper for a detailed exposition of this approach.
 \par 
In Problem P1, until the time $\tau^1$, the information available to O2 consists only of the messages sent from O1. This is the same information that the designer uses to select the thresholds. Thus O2 can be thought of as playing the role of the designer in the proof of Theorem 5. The fact that the problem of choosing the thresholds can be viewed from O2's perspective is crucial in determining the nature of the information state for this problem. The form of our information state and the approach of viewing the problem from O2's perspective imitates the information state and the philosophy adopted in \cite{Walrand83} for a real-time point-to-point communication problem with noiseless feedback, where the problem of choosing the encoding functions can be viewed from the decoder's perspective.
 
\subsection{Sequential Decomposition for Problem P2}
In this section, we present a sequential decomposition similar to Theorem 5 for Problem P2. In Problem P2, both observers start taking measurements at time \(t=1\). Moreover, O2 is allowed to stop before receiving the final message from O1 (see the time-ordering in Fig.2 for t=1,2,...). In Problem P2, the messages sent from O1 to O2 are still common information among the two observers. The problem of choosing the optimal thresholds for the two observers can still be viewed as a sequential problem from the perspective of a designer who at any time $t$ knows the common information. At each time $t$, the designer must specify the thresholds to be used by observer 1 in case the final message has not already been sent. It also has to specify -for each realization of messages from O1- the set of thresholds to be used at O2. In other words, at each time $t$, one can think that the designer knows the messages sent from O1 to O2 and the designer must choose the thresholds to be used by O1 and O2 at time $t$. The designer's problem can therefore be viewed as a sequential optimization problem with classical information structure.
\par
Unlike Problem P1, O2's information no longer coincides with the designer's information of all previous messages from O1, since O2 has its own observations as well. The fact that the designer's problem can no longer be viewed from O2's perspective implies that the information state found for Problem P1 is no longer works for this problem. The main challenge now is to find a suitable information state sufficient for performance evaluation for the designer's problem.  We present such an information state and the resulting dynamic program below. 
 \par
  As mentioned earlier, once observer 1 has sent its final message to observer 2, the optimization problem for observer 2 becomes the well known centralized sequential detection problem studied by Wald. The thresholds characterizing the optimal policy and the cost of the optimal policy are known. For a Wald problem with horizon \(T\) and a starting belief \(\pi\) on the event \(\{H=0\}\), the cost of using the optimal Wald thresholds is a function of the belief \(\pi\) which we denote by \(K^{T}(\pi)\). The designer's task is to select the sequence of thresholds to be used by observer 1 and the sequence of thresholds to be used by observer 2 until the final message has been sent from O1 to O2. After O1's final message has been sent, O2's thresholds are known to be the Wald thresholds with appropriate time-horizon. We will now present a sequential decomposition for the designer.
\par
Recall that we defined observer 2's belief on $H$ as follows:
\[ \pi^2_t(Y^2_{1:t},Z^1_{1:t}) := P^{\Gamma^1}(H=0|Y^2_{1:t},Z^1_{1:t}) \]
$\pi^2_t$ evolves in time as O2 gets more measurements and messages. Once O2 has announced its final decision, its belief on $H$ does not change with time (since O2 is no longer making measurements or listening to messages from O1). 
 We begin with the following definition and lemma.
\begin{definition}
For \(t=1,2,...,T^1\) and a given choice of observer 1 and observer 2's strategies from time instant \(1\) to \(t-1\), (that is, \(\Gamma^1_{t-1} =(\gamma^1_1,\gamma^1_2,...,\gamma^1_{t-1})\) and \(\Gamma^2_{t-1}=(\gamma^2_1,\gamma^2_2,...,\gamma^2_{t-1})\)), define
\[D_t := \mathbbm{1}_{\tau^2 \geq t}\]
 \[ \psi_t := P^{\Gamma^1_{t-1},\Gamma^2_{t-1}}(H,\pi^1_t, \pi^2_{t-1}, D_t|Z^1_{1:t-1}=b_{1:t-1}) \]
 where \(\tau^2\) is the stopping time of O2 as defined in (\ref{eq:Pf5}).
 For \(t=1,2,...,T^1-1\) and for a given choice of strategies (\(\Gamma^1_t=(\gamma^1_1,\gamma^1_2,...,\gamma^1_{t})\)) and (\(\Gamma^2_{t-1}=(\gamma^2_1,\gamma^2_2,...,\gamma^2_{t-1})\)), define
\[ \phi_t[z^1_{t}] := P^{\Gamma^1_{t},\Gamma^2_{t-1}}(H,\pi^1_t, \pi^2_{t}, D_t|Z^1_{1:t-1}=b_{1:t-1},Z^1_{t}=z^1_{t}) \]
where \(z^1_t \in \{0,1,b\}\).
\end{definition}


\vspace{10pt}
\begin{lemma} Consider any policy \(\gamma^1_t\), \(t=1,2,...,T^1\) for observer 1 that is characterized by 4 thresholds \((\alpha^1_t,\beta^1_t,\delta^1_t,\theta^1_t)\), for \(t=1,2,...,T^1-1\) and a threshold \(\alpha^1_{T^1}\) at time \(T^1\) (Theorem 2), and a policy \(\gamma^2_t\),   \(t=1,2,...,T^2\) which is characterized by thresholds \((\alpha^2_t,\beta^2_t)\), \(t=1,2,...,T^1-1\) to be used if O1 has not sent a final message and the Wald thresholds \((w^1_t,w^2_t), t=1,2,....,T^2\) to be used if the final message from O1 has been received. Then, we have:
\begin{enumerate}[i)]
\item There exist transformations \(Q^1_t\)  for \(t=1,2,...,T^1\) such that
\[\phi_t[z^1_t] = Q^1_t(\psi_t, \gamma^1_t,z^1_t)\]
for \(z^1_t \in \{0,1,b\}\), and 
\item There exist transformations \(Q^2_t\), \(t=1,2,...,T^1-1\) such that 
\[\psi_{t+1} = Q^2_t(\phi_t[b],\gamma^2_t)\]
\end{enumerate}
\end{lemma}
\begin{proof}
See Appendix E.
\end{proof}
\par
We can now present a sequential decomposition of problem P2.
\begin{theorem}
 For \(t=1,2,...,T^1-1\), there exist functions \(\mathcal{F}_{t}(\psi_{t},\alpha^1_{t},\beta^1_{t},\delta^1_t,\theta^1_t)\) and \(\mathcal{F}^*_t(\psi_{t})\) and \(\mathcal{G}_t(\phi_{t}[b],\alpha^2_t, \beta^2_t)\) and \(\mathcal{G}^*_t(\phi_{t}[b])\) where
 \[  \mathcal{F}^{*}_{t}(\psi_t) = \inf_{\alpha^1_{t},\beta^1_{t},\delta^1_t,\theta^1_t}\mathcal{F}_{t}(\psi_{t},\alpha^1_{t},\beta^1_{t},\delta^1_t,\theta^1_t)\]
  \[  \mathcal{G}^{*}_{t}(\phi_{t}[b]) = \inf_{\alpha^2_{t},\beta^2_{t}}\mathcal{G}_{t}(\phi_{t}[b],\alpha^2_{t},\beta^2_{t})\] 
 and for \(t=T^1\), there exist functions \(\mathcal{F}(\psi_{T^1},\alpha^1_{T^1})\) and \(\mathcal{F}^{*}_{T^1}(\psi_{T^1})\) where
 \[ \mathcal{F}^{*}_{T^1}(\psi_{T^1}) = \inf_{\alpha^1_{T^1}}\mathcal{F}(\psi_{T^1},\alpha^1_{T^1})\]
 such that the optimal thresholds can be evaluated from these functions as follows:
 \begin{enumerate}
 \item Note that \(\psi_1\) is fixed a priori and does not depend on any design choice. The optimal thresholds at \(t=1\) for O1 are given by optimizing parameters in the definition of \(\mathcal{F}^*_1(\psi_1)\). \\
 \item Once O1's thresholds are fixed, \(\phi_1[b]\) is fixed by Lemma 4. The optimizing thresholds to be used by O2 if a blank message was received are given by optimizing parameters in  the definition of \(\mathcal{G}^*_1(\phi_1[b])\). In case a \(0\) or \(1\) was receiver from O1, the optimal thresholds for O2 from this time onwards are the Wald thresholds for a finite horizon \(T^2-1\).\\
\item Continuing sequentially, \(\psi_t\) is fixed by the choice of past thresholds and the optimal thresholds for O1 at time \(t\) are given  by optimizing parameters in the definition of \(\mathcal{F}^*_t(\psi_t)\). Once O1's thresholds are fixed, \(\phi_t[b]\) is fixed by lemma 4. The optimizing thresholds to be used by O2 if a blank message was received are given by optimizing parameters in  the definition of \(\mathcal{G}^*_t(\phi_t[b])\). In case a \(0\) or \(1\) was receiver from O1, the optimal thresholds for O2 from this time onwards are the Wald thresholds for a finite horizon \(T^2-t\).
\end{enumerate}  
\end{theorem}
\begin{proof}
See Appendix F.
\end{proof} 
\par
As in Theorem 5, the sequential decomposition in Theorem 6 is a dynamic programing result for the designer's sequential problem of choosing the thresholds for O1 and O2.
At time $t$, $\psi_t$ is the designer's information state just before selecting the four thresholds to be used at O1 to decide its message $Z^1_t$, whereas $\phi_t$ is designer's information state just before selecting the thresholds to be used by O2 to decide $U^2_t$. (See Fig. 2). The actual form of the functions $\mathcal{F}_t$ and
$\mathcal{G}_t$ is obtained by backward induction in Appendix F.   
\section{Infinite Horizon Problem}\label{sec:Inf_Horizon}
   In this section we analyze infinite horizon analogues of problems P1 and P2. We first focus on Problem P2.
\subsection{Problem P2 with Infinite Horizon} Consider the model of Problem P2 as described in Section \ref{sec:PF}. We remove the restriction on the boundedness of the stopping times, that is, \(\tau^1\) and \(\tau^2\) need not be bounded. The optimization problem is to select policies \(\Gamma^1 = (\gamma^1_1,\gamma^1_2,...)\) and \(\Gamma^2 = (\gamma^2_1,\gamma^2_2,...)\) to minimize
 \begin{equation}
 \mathds{E}^{\Gamma^1,\Gamma^2}\{c^1\tau^1+c^2\tau^2+J(U^2_{\tau^2},H)\}
 \end{equation}
 where \(\tau^1,\tau^2\) and \(U^2\) are defined by equations (\ref{eq:Pf1}), (\ref{eq:Pf4}) and (\ref{eq:Pf5}). We assume that the cost parameters \(c^1\), \(c^2\) are finite positive numbers and that \(J(U^2,H)\) is non-negative and bounded by a constant \(L\) for all \(U^2\) and \(H\).\\
\emph{Remark:} We can restrict attention to policies for which \(E\{\tau^1\}\) and \(E\{\tau^2\}\) are finite, since otherwise the expected cost would be infinite. Thus, we have that \(\tau^1\) and \(\tau^2\) are almost surely finite. However, the stopping times may not necessarily be bounded even under optimal policies.  
   
\subsubsection{Qualitative Properties for Observer 2} 
Consider any fixed policy \(\Gamma^1\) for Observer 1. We will provide structural results on optimal policies for Observer 2 that hold for any choice of \(\Gamma^1\).  Consider the case when observer 2 has not stopped before time \(t\). Consider a realization of the information available to O2 at time \(t\) - \(y^2_{1:t},z^1_{1:t}\) and let \(\bar{\pi}^2_t=P^{\Gamma^1}(H=0|y^2_{1:t},z^1_{1:t}) \) be the realization of O2's belief on $H$. Let \(\mathcal{A}^{\infty}\) be the set of all policies available to O2 at time \(t\) after having observed \(y^2_{1:t},z^1_{1:t}\), and let  \(\mathcal{A}^{T^2}\) be the subset of policies in \(\mathcal{A}^{\infty}\) for which the stopping time \(\tau^2\) is  less than or equal to a finite horizon \(T^2, (t \leq T^2 < \infty)\). Then, from the analysis for the finite-horizon problem P2, we know that there exist value-functions \(\tilde{V}_t^{T^2}(z^1_{1:t},\bar{\pi}^2_t)\) such that
	\begin{align} &\tilde{V}_t^{T^2}(z^1_{1:t},\bar{\pi}^2_t) \nonumber \\ 
	&= \inf_{\Gamma^2 \in \mathcal{A}^{T^2}} E^{\Gamma^1}[c^1\tau^1 + c^2\tau^2 + J(U^2_{\tau^2},H)|y^2_{1:t},z^1_{1:t}] 
	\end{align} 
	This value-function is the optimal finite horizon cost for observer 2 with horizon $T^2$. \\
	We define the following function:
\begin{align}
 &\tilde{V}_t^{\infty}(z^1_{1:t},y^2_{1:t}) \nonumber \\ 
	&= \inf_{\Gamma^2 \in \mathcal{A}^{\infty}} E^{\Gamma^1}[c^1\tau^1 + c^2\tau^2 + J(U^2_{\tau^2},H)|y^2_{1:t},z^1_{1:t}] 
\end{align} 		

\begin{lemma}
\begin{enumerate}[i)]
\item The value functions \(\tilde{V}_t^{T^2}(z^1_{1:t},\bar{\pi}^2_t)\) are non-increasing in \(T^2\) and bounded below by 0, hence the limit \(\lim_{T^2 \to \infty}\tilde{V}_t^{T^2}(z^1_{1:t},\bar{\pi}^2_t)\) exists.
\item Moreover, 
 \[ \tilde{V}_t^{\infty}(z^1_{1:t},y^2_{1:t}) = \lim_{T^2 \to \infty}\tilde{V}_t^{T^2}(z^1_{1:t},\bar{\pi}^2_t) \]
 \end{enumerate}
\end{lemma}
\begin{proof} See Appendix G.
\end{proof}
We can now prove the following theorem:
\begin{theorem}
 For a fixed policy \(\Gamma^1\) for O1, an optimal policy for O2 is of the form:
 \[ U^2_t = \left \{ \begin{array}{ll}
               1 & \mbox{if $\pi^2_t \leq \alpha_t(Z^1_{1:t})$}  \\
               N & \mbox{if $\alpha_t(Z^1_{1:t}) < \pi^1_t < \beta_t(Z^1_{1:t})$} \\
               0 & \mbox{if $\pi^1_t \geq \beta_t(Z^1_{1:t})$}
               \end{array}
               \right. \] 
               where \(0 \leq \alpha_t(Z^1_{1:t}) \leq \beta_t(Z^1_{1:t}) \leq 1\) are thresholds that depend on the sequence of messages received from O1 (\(Z^1_{1:t}\)).
\end{theorem} 
\begin{proof}
   Consider a realization $y^2_{1:t}$, $z^1_{1:t}$ of O2's observations and messages from O1. Let $\bar{\pi}^2_t$ be the realization of O2's belief, where \(\bar{\pi}^2_t = P^{\Gamma^1}(H=0|z^1_{1:t},y^2_{1:t})\). Since \(\tilde{V}_t^{\infty}(z^1_{1:t},y^2_{1:t}) = \lim_{T^2 \to \infty}\tilde{V}_t^{T^2}(z^1_{1:t},\bar{\pi}^2_t)\), it follows that \(V_t^{\infty}\) is a function only of \(z^1_{1:t}\) and \(\bar{\pi}^2_t\). Since O2 at time \(t\) has only 3 possible choices, we must have: 
   \begin{align}
   &\tilde{V}_t^{\infty}(z^1_{1:t},\bar{\pi}^2_t) := min \{ \nonumber \\
   &E^{\Gamma^1} [J(0,H)|\bar{\pi}^2_t], \nonumber \\
  	                                  &E^{\Gamma^1} [J(1,H)|\bar{\pi}^2_t], \nonumber \\
  	                                  & c^2 + E^{\Gamma^1} [\tilde{V}^{\infty}_{t+1}(Z^1_{1:t+1},\pi^2_{t+1})|\bar{\pi}^2_t,z^1_{1:t}]\} \label{eq:infh3}
  \end{align}
  From Lemma 2, we know that the first two terms are affine in \(\bar{\pi}^2_t\). From Lemma 5, we know that \(\tilde{V}_{t+1}^{\infty}\) is the limit of a sequence of finite-horizon value functions. Now, for a fixed \(z^1_{1:t+1}\), the finite horizon value functions are concave in \(\pi^2_{t+1}\) (from Lemma 2), therefore, for a fixed \(z^1_{1:t+1}\), the limit \(\tilde{V}_{t+1}^{\infty}\) is concave in \(\pi^2_{t+1}\) as well. Using the concavity of \(\tilde{V}_{t+1}^{\infty}\) and following the arguments in the proof of Lemma 2, we can show that the third term in equation (\ref{eq:infh3}) is concave in \(\bar{\pi}^2_t\) for a fixed \(z^1_{1:t}\). Thus, for a given realization of \(z^1_{1:t}\), the infinite horizon value function is minimum of two affine and one concave function. Moreover, it is optimal for O2 to stop if \(\bar{\pi}^2_t =0\) or \(1\). Therefore, the optimal policy for O2 must be of the form:    
  \[ U^2_t = \left \{ \begin{array}{ll}
               1 & \mbox{if $\pi^2_t \leq \alpha_t(Z^1_{1:t})$}  \\
               N & \mbox{if $\alpha_t(Z^1_{1:t}) < \pi^1_t < \beta_t(Z^1_{1:t})$} \\
               0 & \mbox{if $\pi^1_t \geq \beta_t(Z^1_{1:t})$}
               \end{array}
               \right. \] 
\end{proof}
 As in the finite horizon problem, once observer 1 has sent the final message to observer 2, observer 2 is faced with the classical centralized 
 Wald problem. With an infinite horizon, the  optimal Wald policies are characterized by stationary thresholds (say, \((w^1,w^2)\))  that do not change with time \cite{Wald}. Thus, in the infinite horizon version of Problem P2, observer 2's operation can be described by the following algorithm:
 \begin{itemize}
   \item From time \(k=1\) onwards, the optimal policy is to use a threshold rule given by 2 numbers \(\alpha_k(b_{1:k})\) and \(\beta_k(b_{1:k})\), \emph{until O1 sends its final message \(Z^1_k \in \{0,1\}\).} 
   \item From the time O1 sends a final message, start using the stationary Wald thresholds \((w^1,w^2)\).
   \end{itemize} 
\subsubsection{Qualitative Properties for Observer 1}

   Consider a fixed policy \(\Gamma^2\) for O2 which belongs to the set of finite horizon policies \(\mathcal{A}^{T^2}\) with horizon \(T^2\). We will show that given such a policy for O2, Observer 1's infinite horizon optimal policy is characterized by 4 thresholds on its posterior belief. We will employ arguments similar to those used in the previous section.
\par
 Consider the case when observer 1 has not stopped before time \(t\). Consider a realization of the information available to O1 at time \(t\) - \(y^1_{1:t}\) and let $\bar{\pi}^1_t = P(H=0|y^1_{1:t})$ be the realization of O1's belief on $H$. Let \(\mathcal{B}^{\infty}\) be the set of all policies available to O2 at time \(t\) after having observed \(y^1_{1:t}\), and let  \(\mathcal{B}^{T^1}\) be the subset of policies in \(\mathcal{B}^{\infty}\) for which the stopping time \(\tau^1\) is  less than or equal to a finite horizon \(T^1, (t \leq T^1 < \infty)\). Then, from the analysis for the finite-horizon problem P2, we know that there exist value-functions \(V_t^{T^1}(\bar{\pi}^1_t)\) such that
	\begin{align} &V_t^{T^1}(\bar{\pi}^1_t) \nonumber \\ 
	&= \inf_{\Gamma^1 \in \mathcal{B}^{T^1}} E^{\Gamma^2}[c^1\tau^1 + c^2\tau^2 + J(U^2_{\tau^2},H)|y^1_{1:t}] 
	\end{align} 
	where \(\bar{\pi}^1_t = P(H=0|y^1_{1:t})\). \\
	We define the following function:
\begin{align}
 &V_t^{\infty}(y^1_{1:t}) \nonumber \\ 
	&= \inf_{\Gamma^1 \in \mathcal{B}^{\infty}} E^{\Gamma^2}[c^1\tau^1 + c^2\tau^2 + J(U^2_{\tau^2},H)|y^1_{1:t}] 
\end{align} 		

\begin{lemma}
\begin{enumerate}[i)]
\item For a fixed finite-horizon policy of O2, the value functions \(V_t^{T^1}(\bar{\pi}^1_t)\) for O1 are non-increasing in \(T^1\) and bounded below by 0, hence the limit \(\lim_{T^1 \to \infty}V_t^{T^1}(\bar{\pi}^1_t)\) exists.
\item Moreover, 
 \[ V_t^{\infty}(y^1_{1:t}) = \lim_{T^1 \to \infty}V_t^{T^1}(\bar{\pi}^1_t) \]
 \end{enumerate}
\end{lemma}
\begin{proof} See Appendix H.
\end{proof} 
We can now state the following theorem:
\begin{theorem}
 For a fixed finite-horizon policy \(\Gamma^2\) for O2, an optimal policy for O1 is of the form:
 \[ Z^1_t = \left \{ \begin{array}{ll}
               b & \mbox{if $\pi^1_t < \alpha_t$}  \\
               1 & \mbox{if $\alpha_t \leq \pi^1_t \leq \beta_t$} \\
               b & \mbox{if $\beta_t<\pi^1_t < \delta_t$} \\
               0 & \mbox{if $\delta_t \leq \pi^1_t \leq \theta_t$} \\
               b & \mbox{if $\pi^1_t > \theta_t$}
               \end{array}
               \right. \]  
               where \(0 \leq \alpha_t \leq \beta_t \leq \delta_t \leq \theta_t \leq 1\).
\end{theorem} 
\begin{proof}
   Because of the above Lemma, we conclude that \(V_t^{\infty}(y^1_{1:t})\) depends only on the realization \(\bar{\pi}^1_t\) of the belief ($\bar{\pi}^1_t = P(H=0|y^1_{1:t})$). It is, moreover, a concave function of \(\bar{\pi}^1_t\). The result of the theorem follows by using  arguments similar to those in the proof of Lemma 1.
\end{proof}

\begin{theorem}
 There exist globally \(\epsilon\)-optimal policies \(G^1,G^2\) for observers 1 and 2 respectively, such that, \(G^1\) is characterized by 4 time -varying thresholds. 
\end{theorem}
\begin{proof}
 Consider any \(\epsilon/2\)-optimal pair of policies \(\Gamma^1,\Gamma^2\). Then, by arguments used in Lemma 5, we know that there exist a finite horizon policy \(\Gamma^2_{T^2}\) such that the pair \(\Gamma^1,\Gamma^2_{T^2}\) is at most \(\epsilon/2\) worse than \(\Gamma^1,\Gamma^2\). Since \(\Gamma^2_{T^2}\) is a finite horizon policy, by theorem 8, we conclude that O1 can use a 4-threshold rule  without losing any performance with respect to the policies \(\Gamma^1,\Gamma^2_{T^2}\). Thus, we have an \(\epsilon\) optimal pair of policies where O1's policy is characterized by 4 time-varying thresholds.   
\end{proof}

\subsection{Problem P1 with Infinite Horizon}
   The above analysis for infinite horizon version of Problem P2 can be easily specialized to the case of Problem P1. In particular, observer~2's problem is now the classical Wald problem with infinite horizon; thus its optimal policy is characterized by two stationary thresholds. Moreover, the arguments of Lemma 6 and Theorems 8 and 9 can be repeated without any modification to obtain the same qualitative properties for observer 1 in Problem P1. 
\section{Communication with M-ary Alphabet} \label{sec:enhanced_alpha}
   Consider models of Problem P1 or P2 with the following modification: when observer 1 chooses to stop taking measurements and send a message to observer 2, it can choose to send one of \(M\) possible choices from the set: \(\{0,1,...,M-1\}\). Thus, observer 1's  message at time \(t\) to observer 2, which is a function of all its observations,
   \begin{equation}
 Z^1_t=\gamma^1_t(Y^1_{1:t}),
\end{equation}
 belongs to the set \(\{0,1,...,M-1,b\}\), where we use \(b\) for blank message, that is, no transmission. The sequence of functions \(\gamma^1_t,t=1,2,...,\) constitute the \emph{policy} of observer 1. Let \(\tau^1\) be the stopping time when observer sends a final message to observer 2, that is,
\begin{equation} 
 \tau^1=min\{t: Z^1_t \in \{0,1,...,M-1\}\}
\end{equation}
Observer 2's operation and the overall system objective are the same as in problem P1 or P2. Then, we have the following result:
\begin{theorem}
 In Problems P1 or P2 where observer 1 can send one of \(M\) possible final messages, there is no loss of optimality in restricting attention to policies for observer 1 that are of the form:
 \[ Z^1_{T^1} = \left \{ \begin{array}{ll}
               M-1 & \mbox{if $\pi^1_{T^1} \leq \alpha^{M-1}_{T^1}$} \\
               M-2 & \mbox{if $\alpha^{M-1}_{T^1}< \pi^1_{T^1} \leq \alpha^{M-2}_{T^1}$} \\
               ..... \\
               1 & \mbox{if $\alpha^2_{T^1}<\pi^1_{T^1} \leq \alpha^1_{T^1}$} \\
               0 & \mbox{if $\pi^1_{T^1} > \alpha^1_{T^1}$ }
               \end{array}
               \right. \]
               
               where \(0 \leq \alpha^{M-1}_{T^1} \leq \alpha^{M-2}_{T^1} \leq ... \leq \alpha^1_{T^1} \leq 1\) are \(M-1\) thresholds  and for \(k=1,2,..,T^1-1\),
         \[ Z^1_k = \left \{ \begin{array}{ll}
               b & \mbox{if $\pi^1_k < \alpha^{M-1}_k$}  \\
               M-1 & \mbox{if $\alpha^{M-1}_k \leq \pi^1_k \leq \beta^{M-1}_k$} \\
               b & \mbox{if $\beta^{M-1}_k < \pi^1_k < \alpha^{M-2}_k$} \\
               M-2 & \mbox{if $\alpha^{M-2}_k \leq \pi^1_k \leq \beta^{M-2}_k$} \\
               ... \\
               ...\\
               1 & \mbox{if $\alpha^1_k \leq \pi^1_k \leq \beta^1_k$} \\
               b & \mbox{if $\beta^1_k<\pi^1_k < \alpha^0_k$} \\
               0 & \mbox{if $\alpha^0_k \leq \pi^1_k \leq \beta^0_k$} \\
               b & \mbox{if $\pi^1_k > \beta^0_k$}
               \end{array}
               \right. \]  
 where \(0 \leq \alpha^{M-1}_{k} \leq \beta^{M-1}_{k} \leq \alpha^{M-2}_{k} \leq ... \leq \alpha^1_{k} \leq \beta^1_k  \leq \alpha^0_{k} \leq \beta^0_k \leq 1 \) are \(2M\) thresholds.
\end{theorem}
\begin{proof}
  It is straightforward to extend the arguments of Theorem 1 to show  that for a fixed policy of observer 2 optimal policies of observer 1 are functions of its posterior belief \(\pi^1_t\). Similarly, the proof of Lemma 1 can be extended to show that the value function for observer 1 is minimum of \(M\) affine functions of the belief, that represent the expected cost of stopping and sending one of the \(M\) symbols, and 1 concave function of the belief that represents the expected cost of continuing. Taking minimum of affine and concave functions will result in \(M\) intervals of the belief space \([0,1]\) where it is optimal to stop and send one of the \(M\) symbols. If at some time \(t\), the symbols are not ordered in the monotonically decreasing way as specified in the result above, one can permute the symbols  in policies of O1 and O2 at time \(t\) to get the desired ordering without losing performance. 
\end{proof}

\section{Extension to Multiple Sensors} \label{sec:MS}
      In this section, we extend our  results to the case when several peripheral sensors similar to observer 1 in  Problems P1 and P2 are required to send a single final message to a coordinating sensor (similar to O2) which may be taking its own measurements. We show that the peripheral sensors have similar parametric characterizations of their optimal policies as observer O1 in Problems P1 and P2. We obtain a characterization of coordinator's strategy that is similar to that of O2.
      \par  
     Consider a group of N peripheral sensors: S1,S2,...,SN and a coordinating sensor S0. Each sensor can make repeated observations on the random variable \(H\). As before, we assume that conditioned on \(H\), the observations at  different sensors are independent, and the observations made at different time instants  at any sensor are also independent conditioned on \(H\). 
     \begin{figure}[ht]
\begin{center}

\includegraphics[height=5cm,width=7cm]{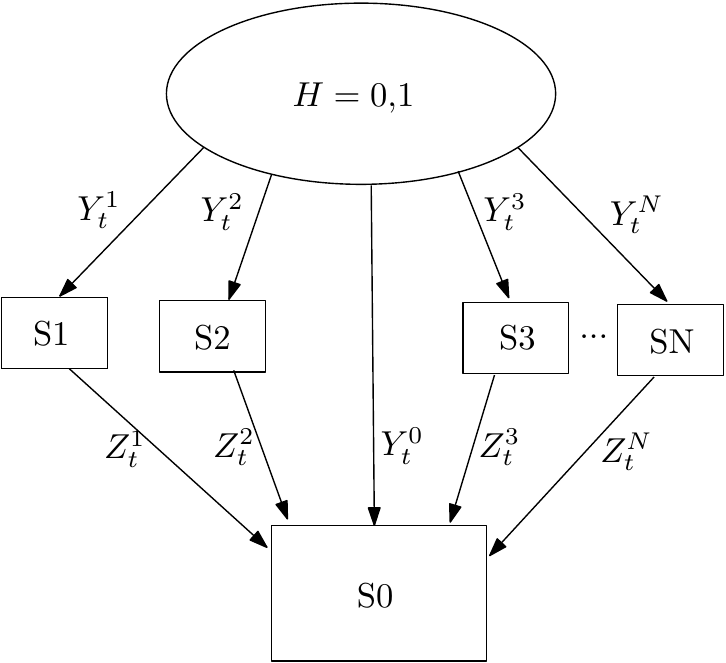}
\caption{Decentralized Detection with N Peripheral Sensors and 1 Coordinating Sensor}

\end{center} 
\end{figure}
     \par
      Each of the peripheral sensors observes its own measurement process \(Y^i_t,i=1,2...,N\) and \(t=1,2,...\). At any time \(t\), the \(i^{th}\) peripheral sensor can decide either to stop and send a binary message \(0\) or \(1\) to the coordinating sensor or to continue taking measurements. Each time the \(i^{th}\) sensor decides to continue taking measurements, a cost \(c^i\) is incurred. Each peripheral sensor sends only a single final message to the coordinator. The policy $\Gamma^i := (\gamma^i_1,\gamma^i_2,\ldots)$ of \(i^{th}\) sensor is of the form:
      \begin{equation}\label{eq:MS1}
      Z^i_t = \gamma^i_t(Y^i_{1:t})
      \end{equation}
      where \(Z^i_t\) is \(i^{th}\) sensor's  message at time \(t\) to the coordinating sensor. \(Z^1_t\) belongs to the set \(\{0,1,b\}\), where we use \(b\) for blank message, that is, no transmission. The time \(\tau^i\) is the stopping time when \(i^{th}\) sensor sends a final message to the coordinating sensor, that is,
\begin{equation} \label{eq:MS1.2}
 \tau^i=min\{t: Z^i_t \in \{0,1\}\}
\end{equation}
The coordinating sensor observes its own measurement process, \(Y^0_t,t=1,2,...\). In addition, it receives messages from all the peripheral sensors (we assume that when the coordinating sensor receives a message it knows which peripheral sensor sent that message). At any time \(t\), S0 can decide to stop and declare a final decision on the hypothesis or take a new measurement and wait for more messages from the peripheral sensors. Each time S0 postpones its decision on the hypothesis, it incurs a cost \(c^0\). When S0 announces a final decision \(U\) on the hypothesis, it incurs a cost given as \(J(U,H)\). Thus, the coordinator's decision at time \(t\) is given as:
\begin{equation}
 U_t = \gamma^0_t(Y^0_{1:t}, Z^1_{1:t},Z^2_{1:t},...,Z^N_{1:t})
\end{equation}
\(U_t\) belongs to the set \(\{0,1,N\}\), where we use \(N\) for a null decision, that is, a decision to continue waiting for more messages and taking more measurements. The sequence of functions \(\Gamma^0 =(\gamma^0_1,\gamma^0_2,...)\) is the policy of the coordinating sensor. The time \(\tau^0\) is the stopping time when S0 announces its final decision on the hypothesis, that is,
\begin{equation} \label{eq:MS5.1}
 \tau^0=min\{t: U_t \in \{0,1\}\}
\end{equation}
We consider the following problem. \\
\emph{Problem P3:} Consider a finite horizon \(T^i\) for the peripheral sensors (that is, we require that \(\tau^i \leq T^i\)) and a finite horizon \(T^0\) for the coordinating sensor, that is, \(\tau^0 \leq T^0\). The optimization problem is to select polices \(\Gamma^0, \Gamma^1,..,.\Gamma^N\) of all the sensors to minimize
\begin{equation}\label{eq:MSobj}
\mathds{E}\{\sum_{i=0}^{N}c^i.\tau^i + J(U_{\tau^0},H)\}
\end{equation}
\par
    We now obtain a characterization of the peripheral sensors' optimal policies. For the \(i^{th}\) peripheral sensor, we define
    \begin{equation}
      \pi^i_t(Y^i_{1:t}) := P(H=0|Y^i_{1:t})
     \end{equation}
\begin{theorem}
      For any peripheral sensor \(i\) and any fixed choice of strategies \(\Gamma^j,\) for \(j=0,1,...,N,j \neq i\), there is an optimal policy of the peripheral sensor \(i\) of the form:
      \[ Z^i_{T^i} = \left \{ \begin{array}{ll}
               1 & \mbox{if $\pi^i_{T^i} \leq \alpha^i_{T^i}$} \\
               0 & \mbox{if $\pi^i_{T^i} > \alpha^i_{T^i}$ }
               \end{array}
               \right. \]
               where \(0 \leq \alpha^i_{T^i} \leq 1\),  and for \(k=1,2,..,T^i-1\),
         \[ Z^i_k = \left \{ \begin{array}{ll}
               b & \mbox{if $\pi^i_k < \alpha^i_k$}  \\
               1 & \mbox{if $\alpha^i_k \leq \pi^i_k \leq \beta^i_k$} \\
               b & \mbox{if $\beta^i_k<\pi^i_k< \delta^i_k$} \\
               0 & \mbox{if $\delta^i_k \leq \pi^i_k \leq \theta^i_k$} \\
               b & \mbox{if $\pi^i_k > \theta^i_k$}
               \end{array}
               \right. \]  
               where \(0 \leq \alpha^i_k \leq \beta^i_k \leq \delta^i_k \leq \theta^i_k \leq 1\). 
\end{theorem}
\begin{proof}
             The main idea of the proof is that once the policies of all sensors except \(i\) are fixed, the optimization problem for the \(i^{th}\) sensor is similar to the problem for O1 in Problem P2. The \(i^{th}\) sensor plays the role of O1 in Problem P2 and the coordinating sensor plays the role of O2. The observations of the coordinating sensor at time \(t\) can be defined as:
             \[ \tilde{Y}^0_t := (Y^0_t,Z^j_t,j=1,2,...,N,j \neq i)\]
             Note that conditioned on \(H\), the observations \(\tilde{Y}^0_t\) are independent of the \(i^{th}\) sensor's observations. We can now follow the arguments of Theorem 1 and 2 to conclude the result for the \(i^{th}\) peripheral sensor. 
\end{proof}
\par
        To find a characterization of the coordinating sensor's policy, we fix the policies of all peripheral sensors and define
        \begin{align} &\pi^0_t(Y^0_{1:t},Z^1_{1:t},Z^2_{1:t},...,Z^N_{1:t}) \nonumber \\ &:= P(H=0|Y^0_{1:t},Z^1_{1:t},Z^2_{1:t},...,Z^N_{1:t}) \end{align}
\begin{theorem}
     For any fixed choice of policies of the peripheral sensors, the policy of the coordinating sensor is given as 
     \[ U_{T^0} = \left \{ \begin{array}{ll}
               1 & \mbox{if $\pi^0_{T^0} \leq \alpha_{T^0}$} \\
               0 & \mbox{if $\pi^0_{T^0} > \alpha_{T^0}$ }
               \end{array}
               \right. \]
   \[ U_k = \left \{ \begin{array}{ll}
               1 & \mbox{if $\pi^0_k \leq \alpha_k(Z^1_{1:k},Z^2_{1:k},...,Z^N_{1:k})$}  \\
               N & \mbox{if $\alpha_k(Z^1_{1:k},Z^2_{1:k},...,Z^N_{1:k}) < \pi^0_k <$} \\ & \mbox{$\beta_k(Z^1_{1:k},Z^2_{1:k},...,Z^N_{1:k})$} \\
               0 & \mbox{if $\pi^0_k \geq \beta_k(Z^1_{1:k},Z^2_{1:k},...,Z^N_{1:k})$}
               \end{array}
               \right. \] 
               where \(0\leq\alpha_k(Z^1_{1:k},Z^2_{1:k},...,Z^N_{1:k})\)\(\leq \beta_k(Z^1_{1:k},Z^2_{1:k},...,Z^N_{1:k})\leq 1\) are thresholds that depend on sequence of messages received from the peripheral sensors.
\end{theorem}
\begin{proof}
  The proof follows the arguments of Theorem 3 and Theorem 4. 
\end{proof}

\section{Conclusion}\label{sec:con}
   We derived structural properties of optimal policies for two observers for a sequential problem in decentralized detection with a single, terminal communication from  observer 1 to the observer 2. It was shown that classical two  threshold rules no longer hold for observer 1.  However, since observer 1's problem is a stopping time problem, a finite parametric characterization of optimal policies is still possible and is described by at most 4 thresholds. A characterization of observer 2's optimal policy was obtained as well. A methodology to find the optimal policies in a sequential manner was presented. We extended the qualitative results to the infinite-horizon versions of the problem, to the problem with increased communication alphabet  and to a related problem with multiple sensors. In all the problems we considered, there is only one message sent from observer 1 to 2. It may still be possible to extend the scope of communication between agents while still satisfying energy and data rate constraints. More general problems where there may be active communication from one observer to the other even before the stopping time remain to be explored.

 \appendices

\section{Proof of Theorem 1}
 Consider an arbitrary choice \(\Gamma^2=(\gamma^2_1,\gamma^2_2,...,\gamma^2_{T2})\) for O2's policy. O2's policy is assumed to be fixed to \(\Gamma^2\) throughout this proof. Note that for a fixed \(\Gamma^2\), \(\tau^2\) and \(U^2_{\tau^2}\) are functions of O2's observation sequence \((Y^2_1,Y^2_2,...,Y^2_{T^2})\) and messages received from O1 \((Z^1_1,...,Z^1_{\tau^1})\). In other words, a policy of O2 induces a stopping time function \(S^{\Gamma^2}\) and an estimate function \(R^{\Gamma^2}\) defined for all possible realizations of the observations of O2 and messages from O1 such that
\begin{equation}
   \tau^2 = S^{\Gamma^2}(Y^2_{1:T^2},Z^1_{1:\tau^1}) \label{eq:Sfunc}
\end{equation}
\begin{equation}
      U^2_{\tau^2} = R^{\Gamma^2}(Y^2_{1:T^2},Z^1_{1:\tau^1}) \label{eq:Rfunc}
\end{equation}
 \par
 Also, by a simple application of Bayes' rule, we know that \(\pi^1_{k+1}\) can be updated from \(\pi^1_k\) and \(Y^1_{k+1}\).
 \begin{align}
   \pi^1_{k+1}= &P(H=0|Y^1_{1:k+1}) \nonumber \\
              = &\frac{P(Y^1_{k+1}|H=0)\pi^1_k}{P(Y^1_{k+1}|H=0)\pi^1_k+P(Y^1_{k+1}|H=1)(1-\pi^1_k)} \label{eq:Ap1}
 \end{align}
 Thus, we have that \begin{equation} \label{eq:belief_update_O1}
  \pi^1_{k+1} = T_k(\pi^1_k,Y^1_{k+1}) 
  \end{equation} where \(T_k\) is defined by (\ref{eq:Ap1}).
 \par
 We will now show that under any policy for O1, the expected future cost at time \(k\) for O1 is lower bounded by the functions \(V_k\) defined in Theorem 1.
 Consider any policy \(\Gamma^1\) for O1. Under the policies \(\Gamma^1\) and \(\Gamma^2\), and for a realization \(y^1_{1:k}\) of O1's observations till time \(k\), let \(W_k(y^1_{1:k})\) be observer 1's expected future cost at time instant \(k\) if it has not sent its final message before time \(k\). That is,
 \begin{align}
 &W_{k}(y^1_{1:k}) :=\mathds{E}^{\Gamma^1,\Gamma^2}[c^1\cdot(\tau^1-k)+c^2\tau^2+J(U^2_{\tau^2},H)|y^1_{1:k}, \nonumber \\&Z^1_{1:k-1}=b_{1:k-1}]
 \end{align}
 First consider time \(T^1\). We have 
 \begin{align}
    &V_{T^1}(\pi) := min\{\nonumber \\
    &\mathds{E}^{\Gamma^2}[c^2\tau^2+J(U^2_{\tau^2},H)|\pi^1_{T^1}=\pi,Z^1_{1:T^1-1}=b_{1:T^1-1},Z^1_{T^1}=0], \nonumber \\
                         &\mathds{E}^{\Gamma^2}[c^2\tau^2+J(U^2_{\tau^2},H)|\pi^1_{T^1}=\pi,Z^1_{1:T^1-1}=b_{1:T^1-1},Z^1_{T^1}=1] \} \label{eq:Apb1}
 \end{align}
  
  If observer 1 has not sent a final decision before time \(T^1\), then under policy \(\Gamma^1\), O1 will either send  \(0\) or \(1\) at time \(T^1\).  O1's expected cost to go at \(T^1\), if it sends a \(0\) at time \(T^1\) is 
 \begin{align}
   &W_{T^1}(y^1_{1:T^1}) = w(y^1_{1:T^1},0) \nonumber \\&:=\mathds{E}^{\Gamma^2}[c^2\tau^2+J(U^2_{\tau^2},H)|y^1_{1:T^1},Z^1_{1:T^1-1}=b_{1:T^1-1},Z^1_{T^1}=0] \label{eq:Ap2} 
 \end{align}
 Similarly, if O1 sends a \(1\) at \(T^1\), its expected cost to go is
  \begin{align}
   &W_{T^1}(y^1_{1:T^1}) = w(y^1_{1:T^1},1) \nonumber \\&:=\mathds{E}^{\Gamma^2}[c^2\tau^2+J(U^2_{\tau^2},H)|y^1_{1:T^1},Z^1_{1:T^1-1}=b_{1:T^1-1},Z^1_{T^1}=1] \label{eq:Ap2.1} 
 \end{align}
 
 Consider the expectation in (\ref{eq:Ap2}). We can write it as 
 \begin{align}
 	&\mathds{E}^{\Gamma^2}[c^2\tau^2+J(U^2_{\tau^2},H)|y^1_{1:T^1},Z^1_{1:T^1-1}=b_{1:T^1-1},Z^1_{T^1}=0] \nonumber \\
 	= &\mathds{E}^{\Gamma^2}[c^2S^{\Gamma^2}(Y^2_{1:T^2},Z^1_{1:\tau^1})+J(R^{\Gamma^2}(Y^2_{1:T^2},Z^1_{1:\tau^1}),H)|y^1_{1:T^1}, \nonumber \\ &Z^1_{1:T^1-1}=b_{1:T^1-1},Z^1_{T^1}=0] \label{eq:apAnew.1}\\
 	=&\mathds{E}^{\Gamma^2}[c^2S^{\Gamma^2}(Y^2_{1:T^2},b_{1:T^1-1},0)+J(R^{\Gamma^2}(Y^2_{1:T^2},b_{1:T^1-1},0),H)\nonumber \\&|y^1_{1:T^1},Z^1_{1:T^1-1}=b_{1:T^1-1},Z^1_{T^1}=0] \label{eq:Ap4}
\end{align}
 where we used (\ref{eq:Sfunc}) and (\ref{eq:Rfunc}) in (\ref{eq:apAnew.1}) and substituted \(Z^1_{1:\tau^1}\) in (\ref{eq:Ap4}) with the values specified in the conditioning term of the expectation. Since the only random variables left in the expectation in (\ref{eq:Ap4}) are \(Y^2_{1:T^2}\) and \(H\), we can write this expectation as
 \begin{align}
 	&\sum\limits_{\substack{h=0,1 \\y^2_{1:T^2} \in \mathcal{Y}^{2}_{1:T^2}}} [ P(y^2_{1:T^2},H=h|\begin{array}{l}y^1_{1:T_1}, Z^1_{1:T^1-1}=b_{1:T^1-1},\\Z^1_{T^1}=0\end{array})\nonumber \\&\times \{c^2S^{\Gamma^2}(y^2_{1:T^2},b_{1:T^1-1},0)+J(R^{\Gamma^2}(y^2_{1:T^2},b_{1:T^1-1},0),h)\}] \label{eq:Ap5}
 \end{align}
 Consider first the term for \(h=0\) in (\ref{eq:Ap5}). Because of the conditional independence of the observations at the two observers, we can write this term as follows:
 \begin{align}
   &\sum\limits_{y^2_{1:T^2} \in \mathcal{Y}^{2}_{1:T^2}} [P(y^2_{1:T^2}|H=0).P(H=0|y^1_{1:T_1}, \nonumber \\&Z^1_{1:T^1-1}=b_{1:T^1-1}, Z^1_{T^1}=0)\times\{c^2S^{\Gamma^2}(y^2_{1:T^2},b_{1:T^1-1},0)+\nonumber \\&J(R^{\Gamma^2}(y^2_{1:T^2},b_{1:T^1-1},0),0)\}] \nonumber \\
   =&\sum\limits_{y^2_{1:T^2} \in \mathcal{Y}^{2}_{1:T^2}} [P(y^2_{1:T^2}|H=0).\pi^1_{T^1}(y^1_{1:T^1})\nonumber\\ &\times \{c^2S^{\Gamma^2}(y^2_{1:T^2},b_{1:T^1-1},0)+J(R^{\Gamma^2}(y^2_{1:T^2},b_{1:T^1-1},0),0)\}] \label{eq:Ap5.1}
 \end{align}
 Similarly, the term for \(h=1\) in (\ref{eq:Ap5}) can be written as,
 \begin{align} \label{eq:Ap6}
  &\sum\limits_{y^2_{1:T^2} \in \mathcal{Y}^{2}_{1:T^2}} [P(y^2_{1:T^2}|H=1).(1-\pi^1_{T^1}(y^1_{1:T^1})) \nonumber \\&\times \{c^2S^{\Gamma^2}(y^2_{1:T^2},b_{1:T^1-1},0)+J(R^{\Gamma^2}(y^2_{1:T^2},b_{1:T^1-1},0),1)\}]
 \end{align}
 Combining equations (\ref{eq:Ap5.1}) and (\ref{eq:Ap6}), we see that the expectation in (\ref{eq:Ap5}) depends on \(\pi^1_{T^1}(y^1_{1:T^1})\) and not on the entire sequence \(y^1_{1:T^1}\). Hence, we can replace \(y^1_{1:T^1}\) by \(\pi^1_{T^1}(y^1_{1:T^1})\) in the conditioning in (\ref{eq:Ap2}). Therefore,  
 \begin{align}
 &w_{T^1}(y^1_{1:T_1},0)\nonumber\\
 = &\mathds{E}^{\Gamma^2}[c^2\tau^2+J(U^2_{\tau^2},H)|y^1_{1:T_1},Z^1_{1:T^1-1}=b_{1:T^1-1},Z^1_{T^1}=0] \nonumber \\
   = &\mathds{E}^{\Gamma^2}[c^2\tau^2+J(U^2_{\tau^2},H)|\pi^1_{T^1}(y^1_{1:T_1}),Z^1_{1:T^1-1}=b_{1:T^1-1}, \nonumber \\&Z^1_{T^1}=0] \nonumber \\
   \geq &V_{T^1}(\pi^1_{T^1}(y^1_{1:T_1})) \label{eq:Ap7}
 \end{align}
 where we used the definition of \(V_{T^1}\) in (\ref{eq:Ap7}). Exactly same arguments can be used if O1 sends a \(1\) at time \(T^1\) to show that
 \begin{align}
 w_{T^1}(y^1_{1:T^1},1):=\mathds{E}^{\Gamma^2}[&c^2\tau^2+J(U^2_{\tau^2},H)|y^1_{1:T^1},\nonumber \\&Z^1_{1:T^1-1}=b_{1:T^1-1}, Z^1_{T^1}=1] \nonumber\\
 &\geq V_{T^1}(\pi^1_{T^1}(y^1_{1:T_1}))
 \end{align}
 Hence, we conclude that the following inequality always holds for policy \(\Gamma^1\):
 \[ W_{T^1}(y^1_{1:T_1}) \geq V_{T^1}(\pi^1_{T^1}(y^1_{1:T^1})) \]

 Now consider time \(k\). Assume that \[ W_{k+1}(y^1_{1:k+1}) \geq V_{k+1}(\pi^1_{k+1}(y^1_{1:k+1})) \] If observer 1 has not sent a final decision before time \(k\), then it will send either a \(0,1\) or \(b\) at time \(k\). Therefore, O1's expected cost to go at \(k\), \(W_k(y^1_{1:k})\), is either
 \begin{align} 
   &w_{k}(y^1_{1:k},0) := \mathds{E}^{\Gamma^2}[c^2\tau^2+J(U^2_{\tau^2},H)|y^1_{1:k},Z^1_{1:k-1}=b_{1:k-1},\nonumber \\&Z^1_{k}=0] \label{eq:Ap7.1}
 \end{align}
 if \(Z^1_k=0\); or
 \begin{align} 
               &w_{k}(y^1_{1:k},1) := \mathds{E}^{\Gamma^2}[c^2\tau^2+J(U^2_{\tau^2},H)|y^1_{1:k},Z^1_{1:k-1}=b_{1:k-1},\nonumber \\&Z^1_{k}=1] \label{eq:Ap8}
 \end{align}
 if \(Z^1_k =1\); or
 \begin{equation}\label{eq:Ap9}
 	  w_{k}(y^1_{1:k},b) := c^1+ \mathds{E}^{\Gamma^2}[W_{k+1}(y^1_{1:k},Y^1_{k+1})|y^1_{1:k},Z^1_{1:k}=b_{1:k}]
 \end{equation}
 if \(Z^1_k= b\). \\
 By arguments similar to those used at time \(T^1\), we can show that \(V_k\) is a lower bound to expressions in (\ref{eq:Ap7.1}) and (\ref{eq:Ap8}). That is,
 \begin{align} \label{eq:Ap7.2}
   &w_{k}(y^1_{1:k},z^1_k)\nonumber \\&=\mathds{E}^{\Gamma^2}[c^2\tau^2+J(U^2_{\tau^2},H)|y^1_{1:k},Z^1_{1:k-1}=b_{1:k-1},Z^1_{k}=z^1_k] \nonumber \\&\geq V_k(\pi^1_k(y^1_{1:k}))
 \end{align}
for \(z^1_k \in \{0,1\}\). \\
Consider equation (\ref{eq:Ap9}). From the induction hypothesis at time \(k+1\), we have that
 \[ W_{k+1} (y^1_{1:k+1}) \geq V_{k+1}(\pi^1_{k+1}(y^1_{1:k+1})) \]
 which implies
 \begin{align}
 &\mathds{E}^{\Gamma^2}[W_{k+1}(y^1_{1:k},Y^1_{k+1})|y^1_{1:k},Z^1_{1:k}=b_{1:k}]  \nonumber \\&\geq \mathds{E}^{\Gamma^2}[V_{k+1}(\pi^1_{k+1}(y^1_{1:k},Y^1_{k+1}))|y^1_{1:k},Z^1_{1:k}=b_{1:k}] \nonumber \\
 &= \mathds{E}^{\Gamma^2}[V_{k+1}(T_k(\pi^1_k(y^1_{1:k}),Y^1_{k+1}))|y^1_{1:k},Z^1_{1:k}=b_{1:k}] \label{eq:Ap10}
 \end{align}
 The above expectation is a function of \(\pi^1_k(y^1_{1:k})\) and the conditional probability:
 \[ P(Y^1_{k+1}|y^1_{1:k},Z^1_{1:k}=b_{1:k}) \]
 which can be expressed as:
 \[P(Y^1_{k+1}|H=0).\pi^1_k(y^1_{1:k}) + P(Y^1_{k+1}|H=1).(1-\pi^1_k(y^1_{1:k})) \]
 
 Thus the expectation in (\ref{eq:Ap10}) depends only on \(\pi^1_k(y^1_{1:k})\) and not the entire sequence \(y^1_{1:k}\); Hence, it can be written as:
 \[\mathds{E}^{\Gamma^2}[V_{k+1}(T_k(\pi^1_k(y^1_{1:k}),Y^1_{k+1}))|\pi^1_k(y^1_{1:k}),Z^1_{1:k}=b]\] 
 Equations (\ref{eq:Ap9}) and (\ref{eq:Ap10}) then imply that
 \begin{align}
 &w_k(y^1_{1:k},b) \nonumber \\
 &= c^1+ \mathds{E}^{\Gamma^2}[W_{k+1}(y^1_{1:k},Y^1_{k+1},Z^1_{k+1})|y^1_{1:k},Z^1_{1:k}=b_{1:k}] \nonumber\\&\geq c^1+ \mathds{E}^{\Gamma^2}[V_{k+1}(T_k(\pi^1_k(y^1_{1:k}),Y^1_{k+1}))|\pi^1_k(y^1_{1:k}),Z^1_{1:k}=b_{1:k}] \nonumber \\
&= c^1+ \mathds{E}^{\Gamma^2}[V_{k+1}(\pi^1_{k+1})|\pi^1_k(y^1_{1:k}),Z^1_{1:k}=b_{1:k}] \nonumber \\
 &\geq V_k(\pi^1_k(y^1_{1:k})) \label{eq:Ap11}
 \end{align}
 where we used the definition of \(V_k\) in (\ref{eq:Ap11}).
 From equations (\ref{eq:Ap7.2}) and (\ref{eq:Ap11}), we conclude that the inequality \(W_k(y^1_{1:k}) \geq V_k(\pi^1_k(y^1_{1:k}))\)  is true. Hence, by induction it holds for all \(k=T^1,T^1-1,...,2,1\). Since \(\Gamma^1\) was arbitrary, we conclude that \(V_k\) are lower bounds on the expected cost to go for O1 under any policy for O1 (with O2's policy fixed at \(\Gamma^2\)).
 \par
  A policy \(\Gamma^{*}\) that always selects the minimizing option in the definition of \(V_k\) for each \(\pi\) will achieve the lower bounds \(V_k\) on \(W_k\) with equality for all \(k\).  Note that the total expected cost of policy \(\Gamma^1\) for O1 is \(c^1+E[W_1(Y^1_1)]\) which is greater that \(c^1+E[V_1(\pi^1_1(y^1_{1}))]\) (since we have shown that \(W_1(y^1_1) \geq V_1(\pi^1_1(y^1_1))\)). Thus, we have that \(\Gamma^{*}\) also achieves the lower bound on total expected cost for any policy. Hence, it is  optimal.
  \par
   Thus, an optimal policy is given by selecting the minimizing option in the definition of \(V_k\) at each \(\pi\). This establishes the dynamic program of Theorem 1 and shows that there is an optimal policy of the form:
 \[Z^1_t = \gamma^{*}_{t}(\pi^{1}_{t})\]

 \section{Proof of Lemma 1}
  Consider the first term in definition of \(V_{T^1}\). 
Using functions \(S^{\Gamma^2}\) and \(R^{\Gamma^2}\) from equations (\ref{eq:Sfunc}) and (\ref{eq:Rfunc}) in first term of (\ref{eq:Apb1}), we get
\begin{align}
 	&\mathds{E}^{\Gamma^2}[c^2\tau^2+J(U^2_{\tau^2},H)|\pi^1_{T^1}=\pi,Z^1_{1:T^1-1}=b_{1:T^1-1},\nonumber \\ &Z^1_{T^1}=0] \nonumber \\
 	= &\mathds{E}^{\Gamma^2}[c^2S^{\Gamma^2}(Y^2_{1:T^2},Z^1_{1:\tau^1}) + J(R^{\Gamma^2}(Y^2_{1:T^2},Z^1_{1:\tau^1}),H)\nonumber \\ &|\pi^1_{T^1}=\pi,Z^1_{1:T^1-1}=b_{1:T^1-1},Z^1_{T^1}=0] \nonumber\\
 	=&\mathds{E}^{\Gamma^2}[c^2S^{\Gamma^2}(Y^2_{1:T^2},b_{1:T^1-1},0) +J(R^{\Gamma^2}(Y^2_{1:T^2},b_{1:T^1-1},0),H)\nonumber \\ &|\pi^1_{T^1}=\pi,Z^1_{1:T^1-1}=b_{1:T^1-1},Z^1_{T^1}=0]\label{eq:Apb2}
\end{align}
where we substituted \(Z^1_{1:T^1}\) in (\ref{eq:Apb2}) with the values specified in the conditioning term of the expectation. Since the only random variables left in the expectation in (\ref{eq:Apb2}) are \(Y^2_{1:T^2}\) and \(H\), we can write this expectation as
 \begin{align}
 	&\sum\limits_{\{h=0,1\}} \sum\limits_{y^2_{1:T^2} \in \mathcal{Y}^{2}_{1:T^2}} [ \nonumber \\ &P(y^2_{1:T^2},H=h|\pi^1_{T^1}=\pi,Z^1_{1:T^1-1}=b,Z^1_{T^1}=0) \nonumber \\ &\{c^2S^{\Gamma^2}(y^2_{1:T^2},b_{1:T^1-1},0) \nonumber \\&+J(R^{\Gamma^2}(y^2_{1:T^2},b_{1:T^1-1},0),h)\}] \label{eq:Apb3}
 \end{align}
 Consider first the term for \(h=0\) in (\ref{eq:Apb3}). Because of the conditional independence of the observations at the two observers, we can write this term as follows:
 \begin{align}
   &\sum\limits_{y^2_{1:T^2} \in \mathcal{Y}^{2}_{1:T^2}} [P(y^2_{1:T^2}|H=0).\pi. \nonumber \\ &\{c^2S^{\Gamma^2}(y^2_{1:T^2},b_{1:T^1-1},0)+J(R^{\Gamma^2}(y^2_{1:T^2},b_{1:T^1-1},0),0)\}] \nonumber \\
   =&\pi \times \nonumber \\ &[\sum\limits_{y^2_{1:T^2} \in \mathcal{Y}^{2}_{1:T^2}}P(y^2_{1:T^2}|H=0).\{c^2S^{\Gamma^2}(y^2_{1:T^2},b_{1:T^1-1},0) \nonumber \\ &+J(R^{\Gamma^2}(y^2_{1:T^2},b_{1:T^1-1},0),0)\}] \label{eq:Apb4} \\
   =&\pi \times A^{\Gamma^2}_{T^1} \label{eq:Apb5}
\end{align}
where \(A^{\Gamma^2}_{T^1}\) is the factor multiplying \(\pi\) in (\ref{eq:Apb4}). Note that this factor  depends only on the choice of O2's policy. Similar arguments for the term corresponding to \(h=1\) in (\ref{eq:Apb3}) show that it can be expressed as
\begin{align}
   &(1-\pi) \times B^{\Gamma^2}_{T^1} \label{eq:Apb6}
\end{align}
Equations (\ref{eq:Apb5}) and (\ref{eq:Apb6}) imply that first term of (\ref{eq:Apb1}) is an affine function of \(\pi\), given as \(A^{\Gamma^2}_{T^1}.\pi + B^{\Gamma^2}_{T^1}.(1-\pi)\). Similar arguments hold for the second term of (\ref{eq:Apb1}). Hence, we have that
\[ V_{T^1}(\pi) := min\{ L^{0}_{T^1}(\pi), L^{1}_{T^1}(\pi)\}\]
Since \(V_{T^1}\) is minimum of two affine functions, it is a concave function of \(\pi\). 
\par
We now proceed inductively. Assume that $V_{k+1}$ is a concave function of $\pi$ and consider \(V_k\),
\begin{align}
  &V_{k}(\pi) := min\{ \nonumber \\&\mathds{E}^{\Gamma^2}[c^2\tau^2+J(U^2_{\tau^2},H)|\pi^1_{k}=\pi,Z^1_{1:k-1}=b_{1:k-1},Z^1_{k}=0], \nonumber \\
                          &\mathds{E}^{\Gamma^2}[c^2\tau^2+J(U^2_{\tau^2},H)|\pi^1_{k}=\pi,Z^1_{1:k-1}=b_{1:k-1},Z^1_{k}=1], \nonumber \\
                          &c^1 + E[V_{k+1}(T_k(\pi^1_k,Y^1_{k+1}))|\pi^1_{k}=\pi,Z^1_{1:k}=b]                 \} \label{eq:Apb7}
 \end{align} 
 Repeating the arguments used for \(V_{T^1}\), it can be shown that first two terms in (\ref{eq:Apb7}) are affine functions of \(\pi\). These are the functions \(L^0_k\) and \(L^1_k\) in Lemma 1. To prove that the third term is concave function of \(\pi\), we use the induction hypothesis that \(V_{k+1}\) is a concave function of \(\pi\). Then, \(V_{k+1}\) can be written as an infimum of affine functions
 \begin{equation} V_{k+1}(\pi) = \inf_i \{\lambda_i\pi + \mu_i \} \label{eq:affine} \end{equation}
 Furthermore, last term in (\ref{eq:Apb7}) can be written as:
 \begin{align}
 	&c^1 + E[V_{k+1}(T_k(\pi,Y^1_{k+1}))|\pi^1_{k}=\pi,Z^1_{1:k}=b_{1:k}] \nonumber \\
 	= &c^1 + \sum\limits_{y^1_{k+1} \in \mathcal{Y}^1} [Pr(y^1_{k+1}|\pi^1_k=\pi,Z^1_{1:k}=b_{1:k}). V_{k+1}(T_k(\pi,y^1_{k+1})) ] \label{eq:Apnew1}
 	\end{align}
 		Using the definition of \(T_k\) from equation (\ref{eq:Ap1}), each term in the above summation can be written as
 	\begin{align}
 	 &Pr(y^1_{k+1}|\pi^1_k=\pi,Z^1_{1:k}=b_{1:k})\nonumber \\&.V_{k+1}\left(\frac{P(y^1_{k+1}|H=0).\pi}{Pr(y^1_{k+1}|H=0).\pi+Pr(y^1_{k+1}|H=1).(1-\pi)}\right) \nonumber \\
 	 &= \{Pr(y^1_{k+1}|H=0).\pi+Pr(y^1_{k+1}|H=1).(1-\pi)\}\nonumber \\&.V_{k+1}\left(\frac{P(y^1_{k+1}|H=0).\pi}{Pr(y^1_{k+1}|H=0).\pi+Pr(y^1_{k+1}|H=1).(1-\pi)}\right) \label{eq:lemma1_new1}
 	 \end{align}
 	 Now using the characterization of \(V_{k+1}\) in terms of the affine functions (from equation \ref{eq:affine}) in the equation (\ref{eq:lemma1_new1}), we obtain
 	 \begin{align}
 	 &\inf_i \{\lambda_i.P(y^1_{k+1}|H=0).\pi + \nonumber \\ &(Pr(y^1_{k+1}|H=0).\pi+Pr(y^1_{k+1}|H=1).(1-\pi)).\mu_i \}
 	\end{align}
 Substituting this expression in (\ref{eq:Apnew1}), we obtain
 \begin{align}	
 	&c^1 + \sum\limits_{y^1_{k+1} \in \mathcal{Y}^1}[\inf_i \{\lambda_i.P(y^1_{k+1}|H=0).\pi + \nonumber \\ &(Pr(y^1_{k+1}|H=0).\pi+Pr(y^1_{k+1}|H=1).(1-\pi)).\mu_i \}] \label{eq:Apb8}
 \end{align} 
 Observe that the expression under the infimum is an affine function of \(\pi\). Hence, taking the infimum over \(i\) gives a concave function of \(\pi\) for each \(y^1_{k+1}\). Since the sum of concave functions is concave, the expression in (\ref{eq:Apb8}) is a concave function of \(\pi\). We will call this function \(G_k(\pi)\). Thus, the value function at time \(k\) given by (\ref{eq:Apb7}) can be expressed as:
  \begin{equation}
  V_{k}(\pi) := min\{ L^{0}_k(\pi), L^{1}_k(\pi), G_k(\pi)\} \label{eq:Apb9}
  \end{equation}
  Since \(V_k\) is minimum of a concave and two affine functions, it itself is a concave function. This completes the argument for the induction step and (\ref{eq:Apb9}) now holds for all \(k=(T^1-1),...,2,1\).

\section{Proof of Theorem 3}
 \begin{proof}
  Let \(\Gamma^1=(\gamma^1_1,\gamma^1_2,...,\gamma^1_{T^1})\) be the fixed policy of O1. By definition of \(\pi^2_{k+1}\), we have
   \begin{align}
    &\pi^2_{k+1}(Y^2_{1:k+1},Z^1_{1:k+1}) := P^{\Gamma^1}(H=0|Y^2_{1:k+1},Z^1_{1:k+1}) \nonumber \\
                                            &=\frac{P(H=0,Y^2_{k+1},Z^1_{k+1}|Y^2_{1:k},Z^1_{1:k})}{\sum_{h=0,1}P(H=h,Y^2_{k+1},Z^1_{k+1}|Y^2_{1:k},Z^1_{1:k})}   \label{eq:Apc1}
   \end{align}
   (although we omit the superscript \(\Gamma^1\) for ease of notation, it should be understood that these probabilities are defined with a fixed \(\Gamma^1\).) \\ 
   Consider the numerator in (\ref{eq:Apc1}). It can be written as:
   \begin{align}
   &P(Y^2_{k+1}|H=0,Y^2_{1:k},Z^1_{1:k+1}).P(Z^1_{k+1}|H=0,Y^2_{1:k},Z^1_{1:k}). \nonumber \\&P(H=0|Y^2_{1:k},Z^1_{1:k}) \nonumber \\
   &=P(Y^2_{k+1}|H=0).P(Z^1_{k+1}|H=0,Z^1_{1:k}).\pi^2_k(Y^2_{1:k},Z^1_{1:k}) \label{eq:Apc2}
   \end{align}
   where we used conditional independence of the observations in (\ref{eq:Apc2}). Under a fixed policy of O1, \(Z^1_k\) s are well-defined random variables and hence the second term in (\ref{eq:Apc2}) is well-defined. Similar expressions can be obtained for the terms in the denominator of (\ref{eq:Apc1}). Thus, we have that \(\pi^2_{k+1}\) is a function of \(\pi^2_k\), \(Y^2_{k+1}\) and \(Z^1_{1:k+1}\). That is,
   \begin{equation}
    \pi^2_{k+1} = \tilde{T}_k(\pi^2_k,Y^2_{k+1},Z^1_{1:k+1}) \label{eq:Apc2.1}
   \end{equation}
   In the statement of Theorem 3, we defined \(\tilde{V}_{T^2}\) as
   \begin{align}
  	\tilde{V}_{T^2}(z^1_{1:T^1},\pi) := min \{ &E^{\Gamma^1} [J(0,H)|\pi^2_{T^2}=\pi], \nonumber \\
  	                                  &E^{\Gamma^1} [J(1,H)|\pi^2_{T^2}=\pi]  \} \label{eq:A22.1}
  \end{align}
    
    If O2 has not declared a final decision on the hypothesis till \(T^2-1\), and selects \(U^2_{T^2} =0\), then his future cost at time \(T^2\) is 
    \begin{align}
    &\tilde{W}_{T^2}(y^2_{1:T^2},z^1_{1:T^1},0) := E^{\Gamma^1}[J(0,H)|y^2_{1:T^2},z^1_{1:T^1}] \nonumber \\ = &\pi^2_{T^2}(y^2_{1:T^2},z^1_{1:T^1}).J(0,0) + (1-\pi^2_{T^2}(y^2_{1:T^2},z^1_{1:T^1})).J(0,1) \nonumber \\
    = &E^{\Gamma^1}[J(0,H)|\pi^2_{T^2}(y^2_{1:T^2},z^1_{1:T^1})]
    \end{align}
    which corresponds to the first term in definition of \(\tilde{V}_{T^2}\) at \(\pi^2_{T^2}(y^2_{1:T^2},z^1_{1:T^1})\).
    A similar expression is true if \(U^2_{T^2} =1\). In either case, we have from the definition of \(\tilde{V}_{T^2}\) that for $u \in \{0,1\}$,
    \begin{align} \tilde{W}_{T^2}(y^2_{1:T^2},z^1_{1:T^1},u) &:= E^{\Gamma^1}[J(u,H)|y^2_{1:T^2},z^1_{1:T^1}] \nonumber \\&\geq \tilde{V}_{T^2}(z^1_{1:T^1},\pi^2_{T^2}(y^2_{1:T^2},z^1_{1:T^1})) \end{align}
   thus, the optimal action at time \(T^2\) is to select the minimizing option in the definition of \(\tilde{V}_{T^2}\) and the optimal future cost is the value of \(\tilde{V}_{T^2}\). 
    \par
   We will employ backward induction on the functions \(\tilde{V}_k\) defined in Theorem 3 to show that they represent the optimal value functions for O2.  Consider time instant \(k\). Assume \(\tilde{V}_{k+1}\) gives the optimal cost to go (future cost) function at time \(k+1\). We have, by definition,
   \begin{align}
   &\tilde{V}_k(z^1_{1:k},\pi) := min \{ \nonumber \\
   &E^{\Gamma^1} [J(0,H)|\pi^2_{k}=\pi], \nonumber \\
  	                                  &E^{\Gamma^1} [J(1,H)|\pi^2_{k}=\pi], \nonumber \\
  	                                  & c^2 + E^{\Gamma^1} [\tilde{V}_{k+1}(Z^1_{1:k+1},\pi^2_{k+1})|\pi^2_{k}=\pi,z^1_{1:k}]\} \label{eq:A22.2}
  \end{align}
  At time \(k\), for a realization $y^2_{1:k},z^1_{1:k}$ of O2's observations and O1's messages, the cost of stopping and declaring a decision on the hypothesis at time \(k\) is either
    \begin{align}
     \tilde{W}_k(y^2_{1:k},z^1_{1:k},0) := E^{\Gamma^1}[J(0,H)|y^2_{1:k},z^1_{1:k}] 
   \end{align}   
   or   
    \begin{align}
     \tilde{W}_k(y^2_{1:k},z^1_{1:k},1) := E^{\Gamma^1}[J(1,H)|y^2_{1:k},z^1_{1:k}] 
   \end{align}   
   By arguments similar to those at time \(T^2\), the above terms are the same as the first two terms of \(\tilde{V}_k(z^1_{1:k},\pi^2_k(y^2_{1:k},z^1_{1:k}))\).
     The cost of continuing at time \(k\) is
     \begin{align}
     &\tilde{W}_k(y^2_{1:k},z^1_{1:k},N) = c^2 + E^{\Gamma^1}[\tilde{V}_{k+1}(Z^1_{1:k+1},\pi^2_{k+1})|y^2_{1:k},z^1_{1:k}] \nonumber\\
                                      &= c^2 + \nonumber \\ &E^{\Gamma^1}[\tilde{V}_{k+1}(Z^1_{1:k+1},\tilde{T}_k(\pi^2_k,Y^2_{k+1},Z^1_{1:k+1}))|y^2_{1:k},z^1_{1:k}] \\
                                      &= c^2 + \nonumber \\
                                     &E^{\Gamma^1}[\tilde{V}_{k+1}(z^1_{1:k},Z^1_{k+1},\tilde{T}_k(\pi^2_k,Y^2_{k+1},z^1_{1:k},Z^1_{k+1}))|y^2_{1:k},z^1_{1:k}]  \label{eq:Apc3}
   \end{align}
   The expectation in (\ref{eq:Apc3}) depends on \(\pi^2_k\), \(z^1_{1:k}\) and \(P^{\Gamma^1}(Y^2_{k+1},Z^1_{k+1}|y^2_{1:k},z^1_{1:k}) \). This probability can be written as:
   \begin{align}
    &P(Y^2_{k+1}|H=0).P(Z^1_{k+1}|H=0,z^1_{1:k}). \pi^2_k + \nonumber \\&P(Y^2_{k+1}|H=1).P(Z^1_{k+1}|H=1,z^1_{1:k}).(1- \pi^2_k) 
    \end{align}
   which depends only on \(z^1_{1:k}\) and \(\pi^2_k\). Thus, the cost of continuing is the same as
   \[ c^2 + E^{\Gamma^1}[\tilde{V}_{k+1}(Z^1_{1:k+1},\pi^2_{k+1})|\pi^2_k(y^2_{1:k},z^1_{1:k}), z^1_{1:k}] \]
   which corresponds to the last term in the definition of \(\tilde{V}_k\). Consequently, the optimal action at time $k$ is to select the minimizing option in definition of \(V_k\) and the value of $V_k$ is the optimal expected cost to go at time \(k\). This completes the proof of the assertion of Theorem 3.
\end{proof}   
   
\section{Proof of Lemma 2}
\begin{proof}
 The result of Lemma 2 for time \(T^2\) follows from the definition of \(\tilde{V}_{T^2}\) since
 \[E^{\Gamma^1} [J(0,H)|\pi^2_{T^2}=\pi] = \pi.J(0,0) + (1-\pi).J(0,1)  \]
 This corresponds to the line \(l^0(\pi)\). Similarly,
 \[E^{\Gamma^1} [J(1,H)|\pi^2_{T^2}=\pi] = \pi.J(1,0) + (1-\pi).J(1,1)  \]
 which corresponds to line \(l^1(\pi)\).
Since, for any realization of \(z^1_{1:T^1}\), \(\tilde{V}_{T^2}\) is minimum of two affine functions of \(\pi\), it is concave in \(\pi\) for each \(z^1_{1:T^1}\).  
\par
Assume now that \(\tilde{V}_{k+1}(z^1_{1:k+1},\pi)\) is concave in \(\pi\) for each \(z^1_{1:k+1}\). Then, we can write \(\tilde{V}_{k+1}\) as:
 \begin{equation} \label{eq:Apinfimum2}
 \tilde{V}_{k+1}(z^1_{1:k+1},\pi) = \inf_i \{ \lambda_i(z^1_{1:k+1}).\pi + \mu_i(z^1_{1:k+1}) \} \end{equation}
 
 where \(\lambda_i(z^1_{1:k+1})\) and \(\mu_i(z^1_{1:k+1}) \) are real numbers that depend on \(z^1_{1:k+1}\).
 Consider the value-function at time \(k\).
 \begin{align}
   \tilde{V}_k(z^1_{1:k},\pi) = min \{&E^{\Gamma^1} [J(0,H)|\pi^2_{k}=\pi], \nonumber \\
  	                                  &E^{\Gamma^1} [J(1,H)|\pi^2_{k}=\pi], \nonumber \\
  	                                  & c^2 + E^{\Gamma^1} [\tilde{V}_{k+1}(Z^1_{1:k+1},\pi^2_{k+1})|\pi^2_{k}=\pi,z^1_{1:k}]\}  \label{eq:Apd1}
  \end{align}
  The first two terms in (\ref{eq:Apd1}) correspond to the affine terms \(l^0\) and \(l^1\). The last term in (\ref{eq:Apd1}) can be written as:
  \begin{align}
  &c^2 + E^{\Gamma^1} [\tilde{V}_{k+1}(Z^1_{1:k+1},\pi^2_{k+1})|\pi^2_{k}=\pi,z^1_{1:k}]\}  \nonumber \\
  = &c^2 + E^{\Gamma^1} [\tilde{V}_{k+1}(Z^1_{1:k+1},\tilde{T}_k(\pi^2_k,Y^2_{k+1},Z^1_{1:k+1}))|\pi^2_{k}=\pi,z^1_{1:k}]\}  \nonumber \\
  = &c^2 + \sum\limits_{y^2_{k+1} \in \mathcal{Y}^2} \sum\limits_{z^1_{k+1} \in \{0,1,b\}} [Pr(y^2_{k+1},z^1_{k+1}|\pi^2_k=\pi,z^1_{1:k}). \nonumber \\ &\tilde{V}_{k+1}(z^1_{1:k+1},\tilde{T}_k(\pi,y^2_{k+1},z^1_{1:k+1})) ]  \label{eq:Apd2}	
  \end{align}
  We now use  the fact that \(\tilde{T}_k(\pi,y^2_{k+1},z^1_{1:k+1})\) is given as 
  \begin{align}
  \frac{P(y^2_{k+1}|H=0).P(z^1_{k+1}|H=0,z^1_{1:k}).\pi}{P(y^2_{k+1},z^1_{k+1}|\pi^2_k=\pi,z^1_{1:k})}
  \end{align}
 (see equations (\ref{eq:Apc1}) and (\ref{eq:Apc2})).
 \par
  Focusing on one term of the summation in (\ref{eq:Apd2}) and using (\ref{eq:Apinfimum2}), we can write it as
  \begin{align}
     &P(y^2_{k+1},z^1_{k+1}|\pi^2_k=\pi,z^1_{1:k}) \times \nonumber \\ &\inf_i\{\lambda_i(z^1_{1:k+1}).\left(\frac{P(y^2_{k+1}|H=0).P(z^1_{k+1}|H=0,z^1_{1:k}).\pi}{P(y^2_{k+1},z^1_{k+1}|\pi^2_k=\pi,z^1_{1:k})}\right) \nonumber \\&+ \mu_i(z^1_{1:k+1}) \} \label{eq:Apd3}
  \end{align}
  
  Note that the expression outside the infimum in (\ref{eq:Apd3}) is the same as the denominator in the term multiplying \(\lambda_i(z^1_{1:k+1})\) in (\ref{eq:Apd3}). The expression (\ref{eq:Apd3}) can now be written as
  \begin{align}
  &\inf_i\{\lambda_i(z^1_{1:k+1}).P(y^2_{k+1}|H=0).P(z^1_{k+1}|H=0,z^1_{1:k}).\pi \nonumber \\
 +&\mu_i(z^1_{1:k+1}). P(y^2_{k+1},z^1_{k+1}|\pi^2_k=\pi,z^1_{1:k}) \label{eq:Ape.1}
 \end{align}
   Expanding the probability multiplying \(\mu_i\), we can write (\ref{eq:Ape.1}) as
   \begin{align}
 &\inf_i\{\lambda_i(z^1_{1:k+1}).P(y^2_{k+1}|H=0).P(z^1_{k+1}|H=0,z^1_{1:k}).\pi \nonumber \\ +&\mu_i(z^1_{1:k+1}).(P(y^2_{k+1}|H=0).P(z^1_{k+1}|H=0,z^1_{1:k}).\pi + \nonumber \\ &P(y^2_{k+1}|H=1).P(z^1_{k+1}|H=1,z^1_{1:k}).(1-\pi)) \}  \label{eq:Apd5}
  \end{align} 
 For the given \(z^1_{1:k+1}\) and \(y^2_{k+1}\), the term in the infimum in (\ref{eq:Apd5}) is affine in \(\pi\). Therefore, the expression in (\ref{eq:Apd5}) is concave in \(\pi\). Thus,  for the given realization of \(z^1_{1:k}\), each term in the summation in (\ref{eq:Apd2}) is concave in \(\pi\). Hence, the sum is concave in \(\pi\) as well. This establishes the structure of \(\tilde{V}_k\) in Lemma 2. To complete the induction argument, we only have to note that since \(\tilde{V}_k\) is the minimum of 2 affine and one concave function of \(\pi\) , it is concave in \(\pi\) (for each \(z^1_{1:k}\)).
 \end{proof}
  
\section{Proof of Lemma 4}
\begin{proof}
 We first prove the second part of the lemma.
\par
By definition, we have
\begin{align}
&\psi_{t+1}(h,\pi^1,\pi^2,1) \nonumber \\ 
&= P(H=h,\pi^1_{t+1}=\pi^1, \pi^2_{t}=\pi^2, D_{t+1}=1|Z^1_{1:t}=b_{1:t}) \nonumber \\
&= P(H=h,T_t(\pi^1_t,Y^1_{t+1})=\pi^1, \pi^2_{t}=\pi^2,\nonumber \\ &  D_{t+1}=1|Z^1_{1:t}=b_{1:t}) \label{eq:oplemma.1}
\end{align}
where we used the fact that O1's belief at time \(t+1\) is a function of its belief at time \(t\) and the observation at time \(t+1\), that is, \(\pi^1_{t+1} = T_t(\pi^1_t,Y^1_{t+1})\) (see Appendix A, (\ref{eq:belief_update_O1}). The right hand side (RHS) of (\ref{eq:oplemma.1}) can further be written as:
\begin{align}
& = \int_{y,\pi'}\mathbbm{1}_{T_t(\pi',y)=\pi^1}.P(H=h,\pi^1_t=\pi',Y^1_{t+1}=y,\pi^2_t=\pi^2, \nonumber \\ &D_{t+1}=1|Z^1_{1:t}=b_{1:t}) \nonumber \\
& = \int_{y,\pi'}\mathbbm{1}_{T_t(\pi',y)=\pi^1}.P(Y^1_{t+1}=y|H=h).P(H=h,\pi^1_t=\pi', \nonumber \\
&\pi^2_t=\pi^2,D_{t+1}=1|Z^1_{1:t}=b_{1:t}) \nonumber \\
& =\int_{y,\pi'}\mathbbm{1}_{T_t(\pi',y)=\pi^1}.P(Y^1_{t+1}=y|H=h).P(H=h,\pi^1_t=\pi',\nonumber\\ &\pi^2_t=\pi^2,D_{t}.\mathbbm{1}_{\alpha^2_t<\pi^2_t<\beta^2_t}=1|Z^1_{1:t}=b_{1:t}) \label{eq:ap_optimal.1}
\end{align}
where we used the fact that if  \(Z^1_{1:t}=b_{1:t}\), then the event \(\{\tau^2 \geq t+1\}\) is same as  \(\{\tau^2 \geq t\} \cap \{\alpha^2_t<\pi^2_t<\beta^2_t\}\) and hence \(D_{t+1} = D_t.\mathbbm{1}_{\alpha^2_t<\pi^2_t<\beta^2_t}\). The RHS of (\ref{eq:ap_optimal.1}) can be written as: 
\begin{align}
& =\int_{y,\pi'}\mathbbm{1}_{T_t(\pi',y)=\pi^1}.P(Y^1_{t+1}=y|H=h).P(H=h,\pi^1_t=\pi', \nonumber \\ &\pi^2_t=\pi^2,D_{t}=1|Z^1_{1:t}=b_{1:t}).\mathbbm{1}_{\alpha^2_t<\pi^2<\beta^2_t} \nonumber\\
&= \int_{y,\pi'}\mathbbm{1}_{T_t(\pi',y)=\pi^1}.P(Y^1_{t+1}=y|H=h).\phi_t[b](h,\pi',\pi^2,1) \nonumber \\&.\mathbbm{1}_{\alpha^2_t<\pi^2<\beta^2_t}\label{eq:ap_optimal.2}
\end{align}
The expression given by (\ref{eq:ap_optimal.2}) depends on \(\phi_t[b]\), the thresholds \(\alpha^2_t, \beta^2_t\) specified by \(\gamma^2_t\) and the observation statistics that are known a priori. Thus $\psi_{t+1}(h,\pi^1,\pi^2,D_{t+1}=1)$ is a function of \(\phi_t[b]\) and $\gamma^2_t$.
 \par
Similarly, \(\psi_{t+1}(h,\pi^1,\pi^2,D_{t+1}=0)\) can be written as:
\begin{align}
&\psi_{t+1}(h,\pi^1,\pi^2,D_{t+1}=0) \nonumber \\
& =\int_{y,\pi'}\mathbbm{1}_{T_t(\pi',y)=\pi^1}.P(Y^1_{t+1}=y|H=h).P(H=h,\pi^1_t=\pi',\nonumber\\ &\pi^2_t=\pi^2,D_{t}.\mathbbm{1}_{\alpha^2_t<\pi^2_t<\beta^2_t}=0|Z^1_{1:t}=b_{1:t}) \nonumber \\
&=\int_{y,\pi'}\mathbbm{1}_{T_t(\pi',y)=\pi^1}.P(Y^1_{t+1}=y|H=h).P(H=h,\pi^1_t=\pi',\nonumber\\ &\pi^2_t=\pi^2,D_{t}=0|Z^1_{1:t}=b_{1:t}) \nonumber \\
&+\int_{y,\pi'}\mathbbm{1}_{T_t(\pi',y)=\pi^1}.P(Y^1_{t+1}=y|H=h).P(H=h,\pi^1_t=\pi',\nonumber\\ &\pi^2_t=\pi^2,D_{t}=1|Z^1_{1:t}=b_{1:t}).\mathbbm{1}_{(\pi^2\leq\alpha^2_t)\cup(\pi^2\geq\beta^2_t)} \nonumber \\
&=\int_{y,\pi'}\mathbbm{1}_{T_t(\pi',y)=\pi^1}.P(Y^1_{t+1}=y|H=h).\phi_t[b](h,\pi',\pi^2,0) \nonumber \\
&+\int_{y,\pi'}\mathbbm{1}_{T_t(\pi',y)=\pi^1}.P(Y^1_{t+1}=y|H=h).\phi_t[b](h,\pi',\pi^2,1)\nonumber\\ &.\mathbbm{1}_{(\pi^2\leq\alpha^2_t)\cup(\pi^2\geq\beta^2_t)} \label{eq:ap_optimal.3}
\end{align}
The RHS of (\ref{eq:ap_optimal.3}) depends only on \(\phi_t[b]\) and the thresholds \(\alpha^2_t, \beta^2_t\) specified by \(\gamma^2_t\). This concludes the proof of the second part of the lemma.
\par 
 For the first part of the lemma, consider
\begin{align}
&\phi_t[b](h,\pi^1,\pi^2,1) \nonumber \\ &= P(H=h,\pi^1_t=\pi^1, \pi^2_{t}=\pi^2, D_t=1|Z^1_{1:t}=b_{1:t})  \label{eq:ap3.2}
\end{align}
To simplify this term, first note that
\begin{align}
    &\pi^2_{t}(y^2_{1:t},z^1_{1:t}) := P(H=0|y^2_{1:t},z^1_{1:t}) \nonumber \\                                            &=\frac{P(H=0,y^2_{t},z^1_{t}|y^2_{1:t-1},z^1_{1:t-1})}{\sum_{h=0,1}P(H=h,y^2_{t},z^1_{t}|y^2_{1:t-1},z^1_{1:t-1})}  \label{eq:ap3.0}  
   \end{align}
The numerator in (\ref{eq:ap3.0}) can be written as:
\begin{align}
   &P(y^2_{t}|H=0,y^2_{1:t-1},z^1_{1:t}).P(z^1_{t}|H=0,y^2_{1:t-1},z^1_{1:t-1}). \nonumber \\&P(H=0|y^2_{1:t-1},z^1_{1:t-1}) \nonumber \\
   &=P(y^2_{t}|H=0).P(z^1_{t}|H=0,z^1_{1:t-1}).\pi^2_{t-1}(y^2_{1:t-1},z^1_{1:t-1}) \label{eq:ap3.2.1}
   \end{align}
   where we used the conditional independence of observations given \(H\). Thus the numerator in (\ref{eq:ap3.0}) can be evaluated from $y^2_t,\pi^2_{t-1}$ and $P(z^1_{t}|H=0,z^1_{1:t-1})$. Similar expression can be obtained for the terms in the denominator of (\ref{eq:ap3.0}).
Therefore, we have
 \begin{equation} \label{eq:ap_optimal.5}
\pi^2_{t}(y^2_{1:t},z^1_{1:t}) = \tilde{T}_{t-1}(\pi^2_{t-1},y^2_{t},P(z^1_{t}|H,z^1_{1:t-1}))
\end{equation}
 For \(z^1_{1:t}=b_{1:t}\), we have
\begin{align}
\pi^2_t(y^2_{1:t},b_{1:t}) &= \tilde{T}_{t-1}(\pi^2_{t-1},y^2_{t},P(Z^1_{t}=b|H,Z^1_{1:t-1}=b_{1:t-1})) \nonumber \\
& = \tilde{T}_{t-1}(\pi^2_{t-1},y^2_{t},P(\pi^1_t \in \mathcal{C}_t|H,Z^1_{1:t-1}=b_{1:t-1})) \label{eq:ap_optimal.4}
\end{align}
where \(\mathcal{C}_t := [0,\alpha^1_t) \cup (\beta^1_t,\delta^1_t) \cup (\theta^1_t,1]\).
The conditional probability in the argument of \(\tilde{T}_{t-1}\) is a function of \(\psi_t\) and the thresholds specified by \(\gamma^1_t\). Thus, when $z^1_{1:t}=b_{1:t}$, 
\begin{equation} \label{eq:ap_optimal.4_new}
 \pi^2_t = \tilde{T}_{t-1}(\pi^2_{t-1},Y^2_{t},\psi_t,\gamma^1_t) 
 \end{equation}
(since the function \(\gamma^1_t\) is completely characterized by a set of thresholds, we use \(\gamma^1_t\) to denote the set of thresholds).

Because of (\ref{eq:ap_optimal.4_new}), the expression in (\ref{eq:ap3.2}) can now be expressed as:
\begin{align}
&P(H=h,\pi^1_{t}=\pi^1,\tilde{T}_{t-1}(\pi^2_{t-1},Y^2_{t},\psi_t,\gamma^1_t)=\pi^2,\nonumber \\ &D_{t}=1|Z^1_{1:t}=b_{1:t}) 
\end{align}
This can further be expressed as
\begin{align}
\frac{\begin{array}{l} P(H=h,\pi^1_t=\pi^1, \tilde{T}_{t-1}(\pi^2_{t-1},Y^2_{t},\psi_t,\gamma^1_t)=\pi^2,\\ D_t=1, Z^1_t=b|Z^1_{1:t-1}=b_{1:t-1}) \end{array}}{P(Z^1_t=b|Z^1_{1:t-1}=b_{1:t-1})} \nonumber \\
\nonumber \\
\frac{\begin{array}{l} P(H=h,\pi^1_t=\pi^1, \tilde{T}_{t-1}(\pi^2_{t-1},Y^2_{t},\psi_t,\gamma^1_t)=\pi^2,\\ D_t=1, Z^1_t=b|Z^1_{1:t-1}=b_{1:t-1}) \end{array}}{P(\pi^1_t \in \mathcal{C}_t|Z^1_{1:t-1}=b_{1:t-1})} \label{eq:oplemma.2}
\end{align}
where \(\mathcal{C}_t := [0,\alpha^1_t) \cup (\beta^1_t,\delta^1_t) \cup (\theta^1_t,1]\).
The denominator is a function of a marginal distribution of \(\psi_t\). To simplify the numerator, first note that
\(\psi_t\) is fixed already by the choice of decision functions till time \(t-1\). The numerator in (\ref{eq:oplemma.2}) can therefore be written as:
\begin{align}
& = \int_{y,\pi'}\mathbbm{1}_{\tilde{T}_{t-1}(\pi',y,\psi_t,\gamma^1_t)=\pi^2}.P(H=h,\pi^1_t=\pi^1,Y^2_{t}=y, \nonumber \\ &\pi^2_{t-1}=\pi',  D_{t}=1,Z^1_t=b|Z^1_{1:t-1}=b_{1:t-1}) \nonumber \\
& = \int_{y,\pi'}\mathbbm{1}_{\tilde{T}_{t-1}(\pi',y,\psi_t,\gamma^1_t)=\pi^2}.P(Y^2_t=y|H=h)\nonumber\\&.P(Z^1_t=b|\pi^1_t=\pi^1). P(H=h,\pi^1_t=\pi^1,\pi^2_{t-1}=\pi',\nonumber \\ &D_{t}=1|Z^1_{1:t-1}=b) \nonumber \\
&=\int_{y,\pi'}\mathbbm{1}_{\tilde{T}_{t-1}(\pi',y,\psi_t,\gamma^1_t)=\pi^2}.P(Y^2_t=y|H=h).\mathbbm{1}_{\pi^1 \in \mathcal{C}_t}.\nonumber \\ &P(H=h,\pi^1_t=\pi^1,\pi^2_{t-1}=\pi', D_{t}=1|Z^1_{1:t-1}=b_{1:t-1}) \nonumber \\
&=\int_{y,\pi'}\mathbbm{1}_{\tilde{T}_{t-1}(\pi',y,\psi_t,\gamma^1_t)=\pi^2}.P(Y^2_t=y|H=h).\mathbbm{1}_{\pi^1 \in \mathcal{C}_t}.\nonumber \\ &\psi_t(h,\pi^1,\pi',1) \label{eq:GO.3} 
\end{align}
The expression in (\ref{eq:GO.3}) is a function of \(\psi_t\) and the thresholds specified by \(\gamma^1_t\). Since (\ref{eq:GO.3}) is equal to the numerator of (\ref{eq:oplemma.2}), it follows from (\ref{eq:oplemma.2}) and (\ref{eq:GO.3}), and the fact that the denomination of (\ref{eq:oplemma.2}) is a marginal distribution of $\psi_t$,
\[ \phi_t[b](h,\pi^1,\pi^2,1) = Q^1_t(\psi_t,\gamma^1_t,b) \]
Similar analysis holds for \(D_t=0\) and also for \(\phi_t[0]\) and \(\phi_t[1]\).
\end{proof}
\section{Proof Of Theorem 6}
\begin{proof}
With the appropriate definitions of the information states $\psi_t$ and $\phi_t$, the proof of Theorem 6 is similar to that of Theorem 5. As in the proof of Theorem 5, we proceed backward in time. 

 Consider first the final horizon for O1: \(T^1\). Assume that the designer has already specified functions \(\gamma^1_1,\gamma^1_2,...,\gamma^1_{T^1-1}\) for O1 and \(\gamma^2_1,\gamma^2_2,...,\gamma^2_{T^1-1}\) for O2.  The designer has to select a function to be used by O1 at time \(T^1\) in case O1's final message has not been already sent. By Theorem 2, this function is characterized by a single threshold \(\alpha^1_{T^1}\). The expected future cost for the designer is the cost of a Wald problem with horizon \(T^2-T^1\), if observer 2 has not already declared its final decision. Thus, the expected cost for the designer is:. 
\begin{align}
&\mathds{E}[\{c^2(\tau^2-T^1)+J(U_{\tau^2},H)\}.\mathbbm{1}_{\tau^2 \geq T^1}|Z^1_{1:T^1-1}=b_{1:T^1-1}] \nonumber\\
&=\mathds{E}[K^{T^2-T^1}(\pi^2_{T^1}).\mathbbm{1}_{\tau^2 \geq T^1}|Z^1_{1:T^1-1}=b_{1:T^1-1}] \nonumber \\
&=\mathds{E}[K^{T^2-T^1}(\pi^2_{T^1}).D_{T^1}|Z^1_{1:T^1-1}=b_{1:T^1-1}] \label{eq:GO.2} 
\end{align}
where \(K^{T^2-T^1}(\cdot)\) is the cost of using the optimal Wald thresholds from \(T^1\) onwards with an available time horizon of \(T^2-T^1\). This cost can be expressed as:
\begin{align}
= &\mathds{E}[K^{T^2-T^1}(\pi^2_{T^1}).D_{T^1}|Z^1_{1:T^1-1}=b_{1:T^1-1},Z^1_{T^1}=1]\nonumber \\& \cdot P(Z^1_{T^1}=1|Z^1_{1:T^1-1}=b_{1:T^1-1}) \nonumber \\ 
\nonumber\\
&+ \mathds{E}[K^{T^2-T^1}(\pi^2_{T^1}).D_{T^1}|Z^1_{1:T^1-1}=b_{1:T^1-1},Z^1_{T^1}=0]\nonumber \\& \cdot P(Z^1_{T^1}=0|Z^1_{1:T^1-1}=b_{1:T^1-1}) \nonumber \\
\nonumber \\
= &\int_{\pi^2}[K^{T^2-T^1}(\pi^2).P(\pi^2_{T^1}=\pi^2,D_{T^1}=1|Z^1_{1:T^1-1}=b_{1:T^1-1},\nonumber \\&Z^1_{T^1}=1)].P(\pi^1_{T^1} \leq \alpha^1_{T^1}|Z^1_{1:T^1-1}=b_{1:T^1-1}) \nonumber \\
+&\int_{\pi^2}[K^{T^2-T^1}(\pi^2).P(\pi^2_{T^1}=\pi^2,D_{T^1}=1|Z^1_{1:T^1-1}=b_{1:T^1-1},\nonumber \\&Z^1_{T^1}=0)].P(\pi^1_{T^1} > \alpha^1_{T^1}|Z^1_{1:T^1-1}=b_{1:T^1-1}) \nonumber \\
& =: \mathcal{L}_{T^1}(\phi_{T^1}[0],\phi_{T^1}[1],\psi_{T^1},\alpha^1_{T^1}) \label{eq:GoP2.2}
\end{align}
where we used the fact that the probabilities in the integrals are marginals of \(\phi_{T^1}[1]\) \(\phi_{T^1}[0]\) respectively and the probabilities multiplying the integrals are marginals of \(\psi_{T^1}\). Using Lemma 4, we can write (\ref{eq:GoP2.2}) as
\begin{align}
& = \mathcal{L}_{T^1}(Q^1_{T^1}(\psi^1_{T^1},\alpha^1_{T^1},0),Q^1_{T^1}(\psi^1_{T^1},\alpha^1_{T^1},1),\psi_{T^1},\alpha^1_{T^1})) \nonumber \\
& =: \mathcal{F}_{T^1}(\psi_{T^1},\alpha^1_{T^1})
\end{align}
Thus the optimization problem for the designer is to select \(\alpha^1_{T^1}\) to minimize \(\mathcal{F}_{T^1}(\psi_{T^1},\alpha^1_{T^1})\). Define
\[ \mathcal{F}^{*}_{T^1}(\psi_{T^1}) = \inf_{\alpha^1_{T^1}}\mathcal{F}_{T^1}(\psi_{T^1},\alpha^1_{T^1})\]
For a given \( \psi_{T^1}\), the function \(\mathcal{F}^{*}_{T^1}\) describes the optimal future cost for the designer and the optimizing \(\alpha^1_{T^1}\) gives the best threshold. \\
Proceeding backwards, assume \(\mathcal{F}^{*}_{t+1}\)  describes the designer's future cost from time $t+1$. We now consider the designer's problem of selecting thresholds \(\alpha^2_{t},\beta^2_{t}\) to be used by O2 if it received all blank messages from O1, that is, \(Z^1_{1:t} =b_{1:t}\). The cost at time $t$ is \(J(0,H)\) if observer 2 stops and declares 0, \(J(1,H)\) if observer 2 declares 1. In case, observer 2 does not make a final decision at this point, a cost of $c^2$ is incurred. The future cost for the designer will be the optimal cost at time \(t+1\) which is given by \(\mathcal{F}^*_{t+1}(\psi_{t+1})\). Thus the expected cost is given as:
\begin{align}
&\mathds{E}[c^1(\tau^1-(t+1))+ \{c^2(\tau^2-t)+J(U_{\tau^2},H)\}\cdot D_t\nonumber\\&|Z^1_{1:t}=b_{1:t}] \nonumber\\
&=\mathds{E}[\{J(1,H).\mathbbm{1}_{\pi^2_{t} \leq \alpha^2_{t}} + J(0,H).\mathbbm{1}_{\pi^2_{t} \geq \beta^2_{t}} + \nonumber \\
&c^2\cdot \mathbbm{1}_{\pi^2_{t} \in [\alpha^2_{t},\beta^2_{t}]}\}.D_t|Z^1_{1:t}=b_{1:t}]+ \mathcal{F}^{*}_{t+1}(\psi_{t+1}) \nonumber\\
\nonumber\\ 
 &= \mathds{E}[\{J(1,H).\mathbbm{1}_{\pi^2_{t} \leq \alpha^2_{t}} + J(0,H).\mathbbm{1}_{\pi^2_{t} \geq \beta^2_{t}} + \nonumber \\
&c^2\cdot \mathbbm{1}_{\pi^2_{t} \in [\alpha^2_{t},\beta^2_{t}]}\}.D_t|Z^1_{1:t}=b_{1:t}]+ \mathcal{F}^{*}_{t+1}(Q^2_t(\phi_t[b],\alpha^2_{t},\beta^2_{t})) \label{eq:Gop2.3} \\
 &=: \mathcal{G}_{t}(\phi_{t}[b],\alpha^2_{t},\beta^2_{t})
 \end{align}  
 where we used the fact that the expectation in (\ref{eq:Gop2.3}) depends on  the thresholds \(\alpha^2_{t},\beta^2_{t}\) , and the conditional belief on \(H\), \(D_t\)  and \(\pi^2_{t}\) given \(Z^1_{1:t}=b^1_{1:t}\)- which is a marginal of \(\phi_{t}[b]\).
 Thus the optimization problem for the designer is to select \(\alpha^2_{t},\beta^2_{t}\) to minimize \(\mathcal{G}_{t}(\phi_{t}[b],\alpha^2_{t},\beta^2_{t})\).Define
\[ \mathcal{G}^{*}_{t}(\phi_{t}[b]) = \inf_{\alpha^2_{t},\beta^2_{t}}\mathcal{G}_{t}(\phi_{t},\alpha^2_{t},\beta^2_{t})\] 
\par
Now consider the designer's problem of selecting thresholds \(\alpha^1_{t},\beta^1_{t},\delta^1_{t},\theta^1_{t}\) to be used by O1 at time $t$. The expected future cost is
$K^{T^2-t}(\pi^2_{t})$ if a final message is sent at time $t$ and if O2 had not already stopped (that is, $D_t=1$). In case a blank message is sent, the designer will need to choose thresholds at time $t$ for O2 and the optimal future cost would be given by  $c^1 + \mathcal{G}^{*}_{t}(\phi_t[b])$. The total expected future cost is therefore,
\begin{align}
&\mathds{E}[c^1(\tau^1-t)+ \{c^2(\tau^2-t)+J(U_{\tau^2},H)\}D_t\nonumber\\&|Z^1_{1:t-1}=b_{1:t-1}] \nonumber\\
&=\mathds{E}[K^{T^2-t}(\pi^2_{t}).D_{t}|Z^1_{t}=0,Z^1_{1:t-1}=b_{1:t-1}]\nonumber \\& \cdot P(Z^1_{t}=0|Z^1_{1:t-1}=b_{1:t-1}) \nonumber \\
&+\mathds{E}[K^{T^2-t}(\pi^2_{t}).D_{t}|Z^1_{t}=1,Z^1_{1:t-1}=b_{1:t-1}]\nonumber \\& \cdot P(Z^1_{t}=1|Z^1_{1:t-1}=b_{1:t-1}) \nonumber \\
 &+ [c^1+\mathcal{G}^{*}_{t+1}(\phi_{t}[b])]\cdot P(Z^1_{t}=b|Z^1_{1:t-1}=b_{1:t-1}) \nonumber\\
 \nonumber\\
 =&\mathds{E}[K^{T^2-t}(\pi^2_{t}).D_{t}|Z^1_{t}=0,Z^1_{1:t-1}=b_{1:t-1}] \nonumber \\&\cdot P(\delta^1_{t}<\pi^1_{t}<\theta^1_{t}|Z^1_{1:t-1}=b_{1:t-1}) \nonumber \\
&+\mathds{E}[K^{T^2-t}(\pi^2_{t}).D_{t}|Z^1_{t}=1,Z^1_{1:t-1}=b_{1:t-1}] \nonumber \\&\cdot P(\alpha^1_{t}<\pi^1_{t}<\beta^1_{t}|Z^1_{1:t-1}=b_{1:t-1}) \nonumber \\
 &+[c^1+ \mathcal{G}^{*}_{t+1}(\phi_{t}[b]).D_{t}]\cdot P(\pi^1_{t} \in \mathcal{C}_t|Z^1_{1:t-1}=b_{1:t-1}) \label{eq:last_app_1}\\
 \nonumber\\
 & =:  \mathcal{L}_{t}(\phi_t[0],\phi_t[1],\phi_t[b],\psi_{t},\alpha^1_{t},\beta^1_{t},\delta^1_{t},\theta^1_{t}) \label{eq:last_app_2}\ \\
 & =: \mathcal{F}_{t}(\psi_{t},\alpha^1_{t},\beta^1_{t},\delta^1_{t},\theta^1_{t}) \label{eq:last_app_3}
\end{align}
where, to write (\ref{eq:last_app_2}), we used the fact that the two expectations in (\ref{eq:last_app_1}) are functions of $\phi_t[0]$ and $\phi_t[1]$ (this can be established using analysis similar to that leading to (\ref{eq:GoP2.2})) and the probabilities multiplying the three terms are marginals of $\psi_t$ . Further, since $\phi_t[0]$ and $\phi_t[1]$ are functions of $\psi_t$, we can write (\ref{eq:last_app_2}) as (\ref{eq:last_app_3}). The analysis for time \(t\) can be inductively repeated for all times. 
\end{proof}
\section{Proof of Lemma 5}
 The first part of the lemma follows directly from the fact that \(\tilde{V}_t^{T^2}\) is defined as infimum over a monotonically increasing sequence of sets \(\mathcal{A}^{T^2}\).
 \par
 We will now prove the second part of the lemma.
 \(\tilde{V}^{\infty}_t\) is defined as infimum of the objective over the set of policies \(\mathcal{A}^{\infty}\) which contains \(\mathcal{A}^{T^2}, \forall T^2\), hence we conclude that 
 \begin{equation} \label{eq:inf_new1}
  \tilde{V}_t^{\infty}(z^1_{1:t},y^2_{1:t}) \leq \lim_{T^2 \to \infty}\tilde{V}_t^{T^2}(z^1_{1:t},\bar{\pi}^2_t)
  \end{equation}
   Assume that the inequality in (\ref{eq:inf_new1}) is strict. Then, there exists a policy \(G \in \mathcal{A}^{\infty}\) for observer 2 such that the expected cost under \(G\),  
\[ W_t(G) := E^{\Gamma^1,G}[c^1\tau^1 + c^2\tau^2 + J(U^2_{\tau^2},H)|y^2_{1:t},z^1_{1:t}],\]
is strictly less than \(\lim_{T^2 \to \infty}\tilde{V}_t^{T^2}(z^1_{1:t},\bar{\pi}^2_t)\). Therefore, the policy \(G\) is better than any finite horizon policy. We will now construct a sequence of finite horizon policies \(G_{T^2}, T^2=t,t+1,t+2,...\) such that the expected cost of \(G_{T^2}\) approaches the expected cost of policy \(G\) as \(T^2 \to \infty\). This will contradict the fact that \(W_t(G) < \lim_{T^2 \to \infty}\tilde{V}_t^{T^2}(z^1_{1:t},\bar{\pi}^2_t)\).
Let \(\tau^G\) and \(U_{\tau^G}\) be the stopping time and the decision at the stopping time induced under policy \(G\). The policy \(G_{T^2}\) is characterized by the stopping time \(\tau'\) and the decision at stopping time \(U_{\tau'}\) it induces as follows:
\[\hspace{20pt} \tau' = \left \{ \begin{array}{ll}
               \tau^G & \mbox{if $\tau^G \leq T^2$} \\
                T^2 & \mbox{if $\tau^G > T^2$}
               \end{array}
               \right. \]
               
and
\[\hspace{20pt} U_{\tau'} = \left \{ \begin{array}{ll}
               U_{\tau^G} & \mbox{if $\tau^G \leq T^2$} \\
                0 & \mbox{if $\tau^G > T^2$}
               \end{array}
               \right. \]
Note that \(G_{T^2}\) is finite horizon policy since it always stops no later than the horizon \(T^2\). Define
\[ W_t(G_{T^2}) := E^{\Gamma^1,G_{T^2}}[c^1\tau^1 + c^2\tau^2 + J(U^2_{\tau^2},H)|y^2_{1:t},z^1_{1:t}]\]
 By assumption, the cost under policy \(G\) is better than the cost under any finite horizon policy. Therefore, \(W_t(G_{T^2}) \geq W_t(G)\). Moreover,
\begin{align}
 &W_t(G_{T^2}) - W_t(G)\nonumber \\ &= E^{\Gamma^1,G_{T^2}}[c^1\tau^1 + c^2\tau^2 + J(U^2_{\tau^2},H)|y^2_{1:t},z^1_{1:t}] \nonumber \\
                       &-E^{\Gamma^1,G}[c^1\tau^1 + c^2\tau^2 + J(U^2_{\tau^2},H)|y^2_{1:t},z^1_{1:t}]\nonumber\\
                       & = E[c^2(\tau'-\tau^G)+J(U_{\tau'},H)-J(U_{\tau^G},H)|y^2_{1:t},z^1_{1:t}] \nonumber \\
                       &= E[\{c^2(\tau'-\tau^G) + \nonumber \\
                       &J(U_{\tau'},H)-J(U_{\tau^G},H)\}.\mathbbm{1}_{\tau^G \leq T^2}|y^2_{1:t},z^1_{1:t}] \nonumber \\
                       &+ E[\{c^2(\tau'-\tau^G) + \nonumber \\
                       &J(U_{\tau'},H)-J(U_{\tau^G},H)\}.\mathbbm{1}_{\tau^G > T^2}|y^2_{1:t},z^1_{1:t}] \label{eq:infh1} 
\end{align}
 The first expectation in equation (\ref{eq:infh1}) is \(0\) since  for \(\tau^G \leq T^2\), the policy \(G_{T^2}\) has the same stopping time and the final decision as policy \(G\). Thus, we get:
\begin{align}
&W_t(G_{T^2}) - W_t(G) \nonumber \\ = &E[\{c^2(\tau'-\tau^G)+J(U_{\tau'},H)-J(U_{\tau^G},H)\}.\mathbbm{1}_{\tau^G > T^2}|y^2_{1:t},z^1_{1:t}] \nonumber \\
                      = &E[\{c^2(T^2-\tau^G)\nonumber \\+ &J(0,H)-J(U_{\tau^G},H)\}.\mathbbm{1}_{\tau^G > T^2}|y^2_{1:t},z^1_{1:t}] \nonumber \\
                      \leq &E[\{J(0,H)-J(U_{\tau^G},H)\}.\mathbbm{1}_{\tau^G > T^2}|y^2_{1:t},z^1_{1:t}] \label{eq:infh2} \\
                      \leq &L.E[\mathbbm{1}_{\tau^G > T^2}|y^2_{1:t},z^1_{1:t}] \nonumber \\ = &L.P(\tau^G > T^2|y^2_{1:t},z^1_{1:t})
\end{align}
,where $L$ is the finite positive constant that upper-bounds $J(U,H)$. Since the stopping time under policy \(G\) is almost surely finite (otherwise cost of policy would be infinite), we have that  \( P(\tau^G > T^2|y^2_{1:t},z^1_{1:t}) \to 0\), as \(T^2 \to \infty\). Thus, for any \(\epsilon >0\), there exists a horizon \(T^2\) large enough for which \(W_t(G_{T^2}) - W_t(G) \leq \epsilon\). Therefore,
\[ \lim_{T^2 \to \infty} W_t(G_{T^2}) =W_t(G) \]
Hence, we conclude that there does not exist any policy \(G \in \mathcal{A}^{\infty}\) for which \(W_t(G) < \lim_{T^2 \to \infty}\tilde{V}_t^{T^2}(z^1_{1:t},\bar{\pi}^2_t)\). Therefore, \(\tilde{V}_t^{\infty}(z^1_{1:t},y^2_{1:t}) = \lim_{T^2 \to \infty}\tilde{V}_t^{T^2}(z^1_{1:t},\bar{\pi}^2_t)\).

\section{Proof of Lemma 6}

 The fist part of the lemma follows directly from the fact that \(V_t^{T^2}\) is defined as infimum over a monotonically increasing sequence of sets \(\mathcal{B}^{T^2}\).
 \par
We will now prove the second part of the lemma. Since \(\mathcal{B}^{\infty}\) contains \(\mathcal{B}^{T^1}, \forall T^1\), we conclude that \( V_t^{\infty}(y^1_{1:t}) \leq \lim_{T^1 \to \infty}V_t^{T^1}(\bar{\pi}^1_t)\). Assume that the inequality is strict. Then, there exists a policy \(\Lambda \in \mathcal{B}^{\infty}\) for observer 1 such that the expected cost under \(\Lambda\),  
\[ W_t(\Lambda) := E^{\Lambda,\Gamma^2}[c^1\tau^1 + c^2\tau^2 + J(U^2_{\tau^2},H)|y^1_{1:t}],\]
is strictly less than \(\lim_{T^1 \to \infty}V_t^{T^1}(\bar{\pi}^1_t)\). Therefore, the policy \(\Lambda\) is better than any finite horizon policy. We will now construct a sequence of finite horizon policies \(\Lambda_{T^1}, T^1=t,t+1,t+2,...\) such that the expected cost of \(\Lambda_{T^1}\) approaches the cost of policy \(\Lambda\) as \(T^1 \to \infty\). This will contradict the fact that \( W_t(\Lambda) < \lim_{T^1 \to \infty}V_t^{T^1}(\bar{\pi}^1_t)\). \\
Let \(\tau^\Lambda\) and \(Z^1_{\tau^\Lambda}\) be the stopping time and the decision at the stopping time induced under policy \(\Lambda\). The policy \(\Lambda_{T^1}\) is characterized by the stopping time \(\tau^*\) and the decision at stopping time \(Z^1_{\tau^*}\) it induces as follows:
\[\hspace{20pt} \tau^{*} = \left \{ \begin{array}{ll}
               \tau^\Lambda & \mbox{if $\tau^\Lambda \leq T^1$} \\
                T^1 & \mbox{if $\tau^\Lambda > T^1$}
               \end{array}
               \right. \]
               
and
\[\hspace{20pt} Z^1_{\tau^{*}} = \left \{ \begin{array}{ll}
               Z^1_{\tau^\Lambda} & \mbox{if $\tau^\Lambda \leq T^1$} \\
                0 & \mbox{if $\tau^\Lambda > T^1$}
               \end{array}
               \right. \]
Note that \(\Lambda_{T^1}\) is finite horizon policy since it always stops no later than the horizon \(T^1\). Define
\[ W_t(\Lambda_{T^1}) := E^{\Lambda_{T^1},\Gamma^2}[c^1\tau^1 + c^2\tau^2 + J(U^2_{\tau^2},H)|y^2_{1:t},z^1_{1:t}]\]
 By assumption, the cost under policy \(\Lambda\) is better than cost under any finite horizon policy. Therefore, \(W_t(\Lambda_{T^1}) \geq W_t(\Lambda)\). Moreover,
\begin{align}
 &W_t(\Lambda_{T^1}) - W_t(\Lambda)\nonumber \\ &= E^{\Lambda_{T^1},\Gamma^2}[c^1\tau^1 + c^2\tau^2 + J(U^2_{\tau^2},H)|y^1_{1:t}] \nonumber \\
                       &-E^{\Lambda,\Gamma^2}[c^1\tau^1 + c^2\tau^2 + J(U^2_{\tau^2},H)|y^1_{1:t}]\nonumber\\
                       & = E[c^1(\tau^{*}-\tau^\Lambda)|y^1_{1:t}]+E^{\Lambda_{T^1},\Gamma^2}[c^2\tau^2 + J(U^2_{\tau^2},H)|y^1_{1:t}] \nonumber \\
                       &-E^{\Lambda,\Gamma^2}[c^2\tau^2 + J(U^2_{\tau^2},H)|y^1_{1:t}] \label{eq:inf_neweq1} \\
                       &\leq E^{\Lambda_{T^1},\Gamma^2}[c^2\tau^2 + J(U^2_{\tau^2},H)|y^1_{1:t}] \nonumber \\&-E^{\Lambda,\Gamma^2}[c^2\tau^2 + J(U^2_{\tau^2},H)|y^1_{1:t}] \label{eq:inf_neweq2}
                       \end{align}
                       where we used the fact that since \(\tau^* \leq \tau^\Lambda\), the first term in (\ref{eq:inf_neweq1}) is less than or equal to $0$.  Further, (\ref{eq:inf_neweq2}) can be written as:
\begin{align}
                       &E^{\Lambda_{T^1},\Gamma^2}[\{c^2\tau^2 + J(U^2_{\tau^2},H)\}.\mathbbm{1}_{\tau^\Lambda \leq T^1}|y^1_{1:t}] \nonumber \\
                       &-E^{\Lambda,\Gamma^2}[\{c^2\tau^2 + J(U^2_{\tau^2},H)\}\mathbbm{1}_{\tau^\Lambda \leq T^1}|y^1_{1:t}] \nonumber \\
                       &+ E^{\Lambda_{T^1},\Gamma^2}[\{c^2\tau^2 + J(U^2_{\tau^2},H)\}.\mathbbm{1}_{\tau^\Lambda > T^1}|y^1_{1:t}] \nonumber \\
                       &-E^{\Lambda,\Gamma^2}[\{c^2\tau^2 + J(U^2_{\tau^2},H)\}\mathbbm{1}_{\tau^\Lambda > T^1}|y^1_{1:t}] \label{eq:infhO1.1}
\end{align}

 For all realizations where \(\tau^\Lambda \leq T^1\), the policy \(\Lambda_{T^1}\) has the same stopping time and the final decision as policy \(\Lambda\) and hence they both will send the same realization of messages to O2 and hence O2's policy $\Gamma^2$ will produce the same realizations of $\tau^2$ and $U^2_{\tau^2}$. This implies that the first two terms in (\ref{eq:infhO1.1}) are equal. Thus, (\ref{eq:infhO1.1}) becomes 
 \begin{align}
 &E^{\Lambda_{T^1},\Gamma^2}[\{c^2\tau^2 + J(U^2_{\tau^2},H)\}.\mathbbm{1}_{\tau^\Lambda > T^1}|y^1_{1:t}] \nonumber \\
                       &-E^{\Lambda,\Gamma^2}[\{c^2\tau^2 + J(U^2_{\tau^2},H)\}\mathbbm{1}_{\tau^\Lambda > T^1}|y^1_{1:t}]\nonumber \\ 
                       &\leq (c^2.T^2+L).E[\mathbbm{1}_{\tau^\Lambda > T^1}|y^1_{1:t}] \nonumber \\&= (c^2.T^2+L).P(\tau^\Lambda > T^1|y^1_{1:t}) \label{eq:infhO1.1.2}
  \end{align}
 where we used the fact that \(\tau^2\) is bounded by $T^2$ under policy \(\Gamma^2\) by assumption. 
Since the stopping time under policy \(\Lambda\) is almost surely finite (otherwise cost of policy would be infinite), we have that  \( P(\tau^\Lambda > T^1|y^1_{1:t}) \to 0\), as \(T^1 \to \infty\). Thus, from equations (\ref{eq:inf_neweq1})-(\ref{eq:infhO1.1.2}), we conclude that for any \(\epsilon >0\), there exists a horizon $T^1$ large enough such that \(W_t(\Lambda_{T^1}) - W_t(\Lambda) \leq \epsilon\). Therefore,
\[ \lim_{T^1 \to \infty} W_t(\Lambda_{T^1}) =W_t(\Lambda) \]
Hence, we conclude that there does not exist any policy \(\Lambda \in \mathcal{B}^{\infty}\) for which \(W_t(\Lambda) < \lim_{T^1 \to \infty}\tilde{V}_t^{T^1}(\bar{\pi}^1_t)\). Therefore, \(V_t^{\infty}(y^1_{1:t}) = \lim_{T^1 \to \infty}V_t^{T^1}(\bar{\pi}^1_t)\).

\section*{Acknowledgments}
  This research was supported in part by NSF Grant CCR-0325571 and NASA Grant NNX06AD47G.

\bibliographystyle{IEEEtran}
\bibliography{myref}

\begin{thebibliography}{10}
\providecommand{\url}[1]{#1}
\csname url@samestyle\endcsname
\providecommand{\newblock}{\relax}
\providecommand{\bibinfo}[2]{#2}
\providecommand{\BIBentrySTDinterwordspacing}{\spaceskip=0pt\relax}
\providecommand{\BIBentryALTinterwordstretchfactor}{4}
\providecommand{\BIBentryALTinterwordspacing}{\spaceskip=\fontdimen2\font plus
\BIBentryALTinterwordstretchfactor\fontdimen3\font minus
  \fontdimen4\font\relax}
\providecommand{\BIBforeignlanguage}[2]{{%
\expandafter\ifx\csname l@#1\endcsname\relax
\typeout{** WARNING: IEEEtran.bst: No hyphenation pattern has been}%
\typeout{** loaded for the language `#1'. Using the pattern for}%
\typeout{** the default language instead.}%
\else
\language=\csname l@#1\endcsname
\fi
#2}}
\providecommand{\BIBdecl}{\relax}
\BIBdecl

\bibitem{Tenney_Detection}
R.~R. Tenney and {N. R. Sandell Jr.}, ``Detection with distributed sensors,''
  \emph{IEEE Trans. Aerospace Electron. Systems}, vol. AES-17, no.~4, pp.
  501--510, July 1981.

\bibitem{Tsitsiklis_survey}
J.~N. Tsitsiklis, ``Decentralized detection,'' in \emph{Advances in Statistical
  Signal Processing}.\hskip 1em plus 0.5em minus 0.4em\relax JAI Press, 1993,
  pp. 297--344.

\bibitem{Varshney}
P.~K. Varshney, \emph{Distributed Detection and Data Fusion}.\hskip 1em plus
  0.5em minus 0.4em\relax Springer, 1997.

\bibitem{Radner_team}
R.~Radner, ``Team decision problems,'' \emph{The Annals of Math. Statistics},
  vol.~33, no.~3, pp. 857--881, Sept. 1962.

\bibitem{Tsitsiklis_large}
J.~N. Tsitsiklis, ``Decentralized detection by a large number of sensors,''
  \emph{Mathematics of Control, Signals and Systems}, vol.~1, no.~2, pp.
  167--182, 1988.

\bibitem{Chamberland04}
J.-F. Chamberland and V.~V. Veeravalli, ``Asymptotic results for decentralized
  detection in power-constrained wireless sensor networks,'' \emph{IEEE Journal
  on Selected Areas in Communication}, vol.~22, no.~6, pp. 1007--1015, Aug.
  2004.

\bibitem{Ahlswede_hypothesis}
R.~Ahlswede and I.~Csiszar, ``Hypothesis testing with communication
  constraints,'' \emph{IEEE Trans. on Info. Theory}, vol. IT-32, no.~4, pp.
  533--543, 1986.

\bibitem{VBP_Detection}
V.~V. Veeravalli, T.~Basar, and H.~Poor, ``Decentralized sequential detection
  with a fusion center performing the sequential test,'' \emph{IEEE Trans.
  Inform. Theory}, vol.~39, pp. 433--442, Mar. 1993.

\bibitem{Dec_Wald}
D.~Teneketzis and Y.~C. Ho, ``The decentralized wald problem,''
  \emph{Information and Computation, 73}, pp. 23--44, 1987.

\bibitem{LaVigna_86}
A.~LaVigna, {A.M. Makowski}, and {J.S. Baras}, ``A continuous-time distributed
  version of the wald's sequential hypothesis testing problem,'' \emph{Lecture
  Notes in Control and Information Sciences}, vol.~83, pp. 533--543, 1986.

\bibitem{Teneket_Quickest_Det}
D.~Teneketzis and P.~Varaiya, ``The decentralized quickest detection problem,''
  \emph{IEEE Trans. on Automatic Control}, vol. AC-29, no.~7, pp. 641--644,
  July 1984.

\bibitem{Wald}
A.~Wald, \emph{Sequential Analysis}.\hskip 1em plus 0.5em minus 0.4em\relax
  Wiley, New York, 1947.

\bibitem{Ho1980}
{Y. C. Ho}, ``Team decision theory and information structures,'' in
  \emph{Proceedings of the IEEE}, vol.~68, no.~6, 1980, pp. 644--654.

\bibitem{Tsitsiklis_86}
J.~N. Tsitsiklis, ``On threshold rules in decentralized detection,'' in
  \emph{Proceedings of 25th IEEE Conference of Decision and Control}, Dec.
  1986, pp. 232--236.

\bibitem{Witsenhausen-separation}
H.~S. Witsenhausen, ``Separation of estimation and control for discrete time
  systems,'' \emph{Proceedings of the IEEE}, vol.~59, no.~11, pp. 1557--1566,
  Nov. 1971.

\bibitem{NayyarTeneketzis:2008}
A.~Nayyar and D.~Teneketzis, ``On the structure of real-time encoders and
  decoders in a multi-terminal communication system,'' \emph{IEEE Trans. Info.
  Theory}, submitted.

\bibitem{Walrand83}
J.~Walrand and P.~Varaiya, ``Optimal causal coding-decoding problems,''
  \emph{IEEE Trans. Inf. theory}, vol. IT-29, no.~6, pp. 814--820, Nov. 1983.

\end{thebibliography}
\end{document}